\documentclass[a4paper,10pt]{amsart}
\usepackage{amsmath,amsthm,amssymb}

\date{First version: October 13, 2009 - Revised version: January 26, 2010}

\newcommand{\Dom}{\operatornamewithlimits{Dom}}
\newcommand{\iden}{\operatorname{Id}}
\newcommand{\dive}{\operatorname{div}}

\newcommand{\norm}[1]{\left\Vert#1\right\Vert}
\newcommand{\abs}[1]{\left\vert#1\right\vert}
\newcommand{\set}[1]{\left\{#1\right\}}
\newcommand{\Real}{\mathbb{R}}

\newcommand{\PV}{\operatorname{P.V.}}

\newtheorem{thm}{Theorem}[section]
\newtheorem{prop}[thm]{Proposition}
\newtheorem{cor}[thm]{Corollary}
\newtheorem{lem}[thm]{Lemma}
\newtheorem{exams}[thm]{Examples}
\theoremstyle{definition}

\newtheorem{rem}[thm]{Remark}

\numberwithin{equation}{section}

\author[P. R. Stinga]{Pablo Ra\'ul Stinga}
\address{Departamento de Matem\'aticas \\
          Facultad de Ciencias \\
          Universidad Au\-t\'o\-no\-ma de Madrid \\
          28049 Madrid, Spain}
\email{pablo.stinga@uam.es}

\author[J. L. Torrea]{Jos\'e Luis Torrea}
\address{Departamento de Matem\'aticas \\
          Facultad de Ciencias \\
          Universidad Au\-t\'o\-no\-ma de Madrid \\
          28049 Madrid, Spain}
\email{joseluis.torrea@uam.es}

\thanks{Research supported by Ministerio de Ciencia e Innovaci\'{o}n de Espa\~{n}a MTM2008-06621-C02-01}

\keywords{Fractional Laplacian; harmonic oscillator; Harnack's inequality; degenerate Schr\"{o}dinger equation; heat semigroup}

\subjclass[2000]{26A33, 35J10, 35B05, 35J70, 35K05}

\begin{document}

\title[Extension problem and Harnack's inequality]{Extension problem and Harnack's inequality for some fractional operators}

\begin{abstract}
The fractional Laplacian can be obtained as a Dirichlet-to-Neumann map via an extension problem to the upper half space. In this paper we prove the same type of characterization for the fractional powers of second order partial differential operators in some class. We also get a Poisson formula and a system of Cauchy-Riemann equations for the extension. The method is applied to the fractional harmonic oscillator $H^\sigma=(-\Delta+\abs{x}^2)^\sigma$ to deduce a Harnack's inequality. A pointwise formula for $H^\sigma f(x)$ and some maximum and comparison principles are derived.
\end{abstract}

\maketitle


\section{Introduction}

In the last years there has been a growing interest in the study of nonlinear problems involving fractional powers of the Laplace operator $(-\Delta)^\sigma$, $0<\sigma<1$. The fractional Laplacian of a function $f:\Real^n\to\Real$ is defined via Fourier transform as
\begin{equation}\label{Frac Lap Fourier}
\widehat{(-\Delta)^\sigma f}(\xi)=\abs{\xi}^{2\sigma}\widehat{f}(\xi),
\end{equation}
and it can be expressed by the pointwise formula
\begin{equation}\label{Frac Lap Point}
(-\Delta)^\sigma f(x)=c_{n,\sigma}\PV\int_{\Real^n}\frac{f(x)-f(z)}{\abs{x-z}^{n+2\sigma}}~dz,
\end{equation}
where $c_{n,\sigma}$ is a positive constant. Observe from \eqref{Frac Lap Point} that the fractional Laplacian is a nonlocal operator. This fact does not allow to apply local PDE techniques to treat nonlinear problems for $(-\Delta)^\sigma$. To overcome this difficulty, L. Caffarelli and L. Silvestre showed in \cite{Caffarelli-Silvestre CPDE} that any fractional power of the Laplacian can be determined as an operator that maps a Dirichlet boundary condition to a Neumann-type condition via an extension problem. To be more precise, consider the function $u=u(x,y):\Real^n\times[0,\infty)\to\Real$ that solves the boundary value problem
\begin{align}
\label{Ext bound val Lap} u(x,0)&=f(x),&x&\in\Real^n, \\
\label{Ext equation Lap} \Delta_x u+\frac{1-2\sigma}{y}~u_y+u_{yy}&=0,&x&\in\Real^n,~y>0.
\end{align}
Then, up to a multiplicative constant depending only on $\sigma$,
$$-\lim_{y\to0^+}y^{1-2\sigma}u_y(x,y)=(-\Delta)^\sigma f(x).$$
This characterization of $(-\Delta)^\sigma f$ via the local (degenerate) PDE \eqref{Ext equation Lap} was used for the first time in \cite{Caffarelli-Salsa-Silvestre} to get regularity estimates for the obstacle problem for the fractional Laplacian.

To solve \eqref{Ext bound val Lap}-\eqref{Ext equation Lap}, Caffarelli and Silvestre noted that \eqref{Ext equation Lap} can be though as the harmonic extension of $f$ in $2-2\sigma$ dimensions more (see \cite{Caffarelli-Silvestre CPDE}). From there, they established the fundamental solution and, using a conjugate equation, a Poisson formula for $u$. Furthermore, taking advantage of the general theory of degenerate elliptic equations developed by Fabes, Jerison, Kenig and Serapioni in 1982-83, they proved Harnack's estimates for $u$ (and thus for $f$).

Let $\Omega$ be an open subset of $\Real^n$, $n\geq1$, and let $d\eta$ be a positive measure defined on $\Omega$. Consider a linear second order partial differential operator $L$, that we assume to be nonnegative, densely defined, and self-adjoint in $L^2(\Omega,d\eta)$. The fractional powers $L^\sigma$, $0<\sigma<1$, can be defined in a spectral way, see Section \ref{Section Extension}.

The aim of this paper is to describe any fractional power $L^\sigma$ as an operator that maps a Dirichlet condition to a Neumann-type condition via an extension problem as in \cite{Caffarelli-Silvestre CPDE}, developing also the corresponding properties (Poisson formula, fundamental solution, conjugate equation, Cauchy-Riemann equations). With this characterization, the interior Harnack's inequality for any fractional power of one of the most basic Schr\"{o}dinger operators, the harmonic oscillator $H=-\Delta+\abs{x}^2$, is consequently deduced. Besides, we find an explicit pointwise expression for the nonlocal operator $H^\sigma$ that will allow us to get some maximum and comparison principles.

Fractional operators appear in physics, when considering fractional kinetics and anomalous transport \cite{Zaslavsky}.

The heat-diffusion semigroup $\set{e^{-tL}}_{t\geq0}$ generated by $L$ will play a crucial role in our work.

Our first main result is the following.

\begin{thm}\label{Thm:Extension}
Let $f\in\Dom(L^\sigma)$. A solution of the extension problem
\begin{align}
\label{Ext bound value} u(x,0)&=f(x),&\hbox{on }&\Omega;\\
\label{Ext equation} -L_xu+\frac{1-2\sigma}{y}~u_y+u_{yy}&=0,&\hbox{in }&\Omega\times(0,\infty);
\end{align}
is given by
\begin{equation}\label{Extension}
u(x,y)=\frac{1}{\Gamma(\sigma)}\int_0^\infty e^{-tL}(L^\sigma f)(x)e^{-\frac{y^2}{4t}}~\frac{dt}{t^{1-\sigma}},
\end{equation}
and
\begin{equation}\label{Neumann condtn}
\lim_{y\to0^+}\frac{u(x,y)-u(x,0)}{y^{2\sigma}}=\frac{\Gamma(-\sigma)}{4^\sigma\Gamma(\sigma)}L^\sigma f(x)=\frac{1}{2\sigma}\lim_{y\to0^+}y^{1-2\sigma}u_y(x,y).
\end{equation}
Moreover, the following Poisson formula for $u$ holds:
\begin{equation}\label{Poisson}
u(x,y)=\frac{y^{2\sigma}}{4^\sigma\Gamma(\sigma)}\int_0^\infty e^{-tL}f(x)e^{-\frac{y^2}{4t}}~\frac{dt}{t^{1+\sigma}}=\frac{1}{\Gamma(\sigma)}\int_0^\infty e^{-\frac{y^2}{4r}L}f(x)e^{-r}~\frac{dr}{r^{1-\sigma}}.
\end{equation}
\end{thm}

All identities in Theorem \ref{Thm:Extension} are understood in $L^2(\Omega,d\eta)$. Note that a solution $u$ to the degenerate boundary value problem \eqref{Ext bound value}-\eqref{Ext equation} is written explicitly in terms of the heat semigroup $e^{-tL}$ acting on $L^\sigma f$. From here, the Poisson formula \eqref{Poisson} can be immediately obtained (see the proof in Section \ref{Section Extension}), where no fractional power of $L$ is involved. When $L=-\Delta$, the extension result of \cite{Caffarelli-Silvestre CPDE} is recovered (see Examples \ref{Examples}). More properties concerning the Poisson formula are contained in Theorem \ref{Thm:Poisson}. Moreover, \eqref{Poisson} can be derived as in \cite{Caffarelli-Silvestre CPDE} (Remark \ref{Rem:Sol fund y Poisson}): use the fundamental solution (that involves the kernel of the heat semigroup generated by $L$) and an appropriate conjugate equation \eqref{Conj L} to infer the Poisson kernel (see \eqref{Poisson y Fund sol}). The conjugate equation will be studied in detail by defining Cauchy-Riemann equations \eqref{CR Extension} adapted to equation \eqref{Ext equation}. See Section \ref{Section Extension}.

If $L$ has discrete spectrum, i.e. $L\phi_k=\lambda_k\phi_k$, $\lambda_k\geq0$, and $\set{\phi_k}_{k\in\mathbb{N}_0}$ is an orthonormal basis of $L^2(\Omega,d\eta)$, the definition of the fractional power $L^\sigma$ is given in the natural way: if $f\in L^2(\Omega,d\eta)$ has the property that $\sum_k\lambda_k^{2\sigma}\abs{\langle f,\phi_k\rangle}^2=\sum_k\lambda_k^{2\sigma}\abs{\int_\Omega f\phi_k~d\eta}^2<\infty$, then
\begin{equation}\label{en L2}
L^\sigma f=\sum_k\lambda_k^\sigma\langle f,\phi_k\rangle\phi_k,\qquad\hbox{ sum in }L^2(\Omega).
\end{equation}
In Section \ref{Section Unicidad} it is shown that, under this assumption, \eqref{Ext bound value}-\eqref{Ext equation} has a unique solution $u$ (vanishing as $y\to\infty$) such that \eqref{Neumann condtn} holds in the $L^2(\Omega)$-sense. The proof is elementary using orthogonal expansions: just write $u(x,y)=\sum_kc_k(y)\phi_k(x)$ and observe that the coefficients $c_k$ satisfy a Bessel equation. Hence, for the existence and uniqueness in this case, the general theory of degenerate PDE's mentioned above is not needed. This method also gives us local Neumann solutions (see Subsection \ref{Subsection Neumann}).

Let us now turn to the case of the fractional harmonic oscillator. We will be able to define $H^\sigma f$ for all tempered distributions $f$. If $f$ is a function that has also some local regularity then the extension result is true in the classical sense (Theorem \ref{Thm:Extension smooth} and Remark \ref{util}). This last fact is an essential ingredient for the second main result of this article: the interior Harnack's inequality for $H^\sigma$.

\begin{thm}\label{Thm:Harnack}
Let $x_0\in\Real^n$ and $R>0$. Then there exists a positive constant $C$ depending only on $n$, $\sigma$, $x_0$ and $R$ such that
$$\sup_{B_{R/2}(x_0)}f\leq C\inf_{B_{R/2}(x_0)}f,$$
for all nonnegative functions $f:\Real^n\to\Real$ that are $C^2$ in $B_R(x_0)$ and such that $H^\sigma f(x)=0$ for all $x\in B_R(x_0)$.
\end{thm}

The Harnack's inequality is valid for $0<\sigma<1$ and the proof given in Section \ref{Section Harnack} is based (as we already remarked) on Theorem \ref{Thm:Extension smooth} and the Harnack's inequality for degenerate Schr\"{o}dinger operators proved by C. E. Guti\'{e}rrez in \cite{Gutierrez} (this idea is contained in \cite{Caffarelli-Silvestre CPDE} for the case of the fractional Laplacian). The Harnack's inequality for $H$ ($\sigma=1$) follows from general results (see \cite{Trudinger}).

The final part of the paper is devoted to the study of the pointwise expression of the fractional harmonic oscillator and some of its consequences. To that end we collect some previous facts about the fractional Laplacian $(-\Delta)^\sigma$. The natural way to arrive to \eqref{Frac Lap Point} starting from \eqref{Frac Lap Fourier} would be by taking the inverse Fourier transform. However, this path can be avoided if we consider the classical formula for $L^\sigma$ that involves the heat-diffusion semigroup generated by $L$:
\begin{equation}\label{eq:L sigma}
L^\sigma f(x)=\frac{1}{\Gamma(-\sigma)}\int_0^\infty\left(e^{-tL}f(x)-f(x)\right)~\frac{dt}{t^{1+\sigma}}.
\end{equation}
Note that \eqref{eq:L sigma} is motivated by the identity $\lambda^\sigma=\frac{1}{\Gamma(-\sigma)}\int_0^\infty(e^{-t\lambda}-1)\frac{dt}{t^{1+\sigma}}$, $\lambda>0$. When $L=-\Delta$ and $f\in\mathcal{S}$ in \eqref{eq:L sigma}, the Fourier transform recovers \eqref{Frac Lap Fourier}. Furthermore, the formula allows us to obtain (in a very simple way) expression \eqref{Frac Lap Point} with the constant $c_{n,\sigma}$ computed explicitly and in particular to see (Proposition \ref{Prop:Lap sig a 1}) that if a function $f$ is $C^2$ around some $x\in\Real^n$ then
$$\lim_{\sigma\to1^-}(-\Delta)^\sigma f(x)=-\Delta f(x).$$
In Section \ref{Section Hermite} we put $L=H$ in \eqref{eq:L sigma} to get a pointwise formula for $H^\sigma f(x)$ (see Theorem \ref{Thm:H sig point}) and, from there, some maximum and comparison principles for $H^\sigma$.

Throughout this paper $\mathcal{S}$ is the Schwartz class of rapidly decreasing $C^\infty(\Real^n)$ functions, the letter $C$ denotes a constant that may change in each occurrence and it will depend on the parameters involved (whenever it is necessary we point out this dependence with subscripts) and $\Gamma$ stands for the Gamma function. We restrict our attention to $0<\sigma<1$ and, in this range, $\Gamma(-\sigma):=\frac{\Gamma(1-\sigma)}{-\sigma}<0$.

\section{The extension problem}\label{Section Extension}

We begin with the basics of the spectral analysis that will be used throughout this Section. The complete details can be found in \cite[Ch.~12~and~13]{Rudin}. Since $L$ is a nonnegative, densely defined and self-adjoint operator on $L^2(\Omega,d\eta)=L^2(\Omega)$, there is a unique resolution $E$ of the identity, supported on the spectrum of $L$ (which is a subset of $[0,\infty)$), such that
$$L=\int_0^\infty\lambda~dE(\lambda).$$
The identity above is a shorthand notation that means
$$\langle Lf,g\rangle_{L^2(\Omega)}=\int_0^\infty\lambda~dE_{f,g}(\lambda),\qquad f\in\Dom(L),~g\in L^2(\Omega),$$
where $dE_{f,g}(\lambda)$ is a regular Borel complex measure of bounded variation concentrated on the spectrum of $L$, with $d\abs{E_{f,g}}(0,\infty)\leq\norm{f}_{L^2(\Omega)}\norm{g}_{L^2(\Omega)}$. If $\phi(\lambda)$ is a real measurable function defined on $[0,\infty)$, then the operator $\phi(L)$ is given formally by
\begin{equation}\label{phi(L)}
\phi(L)=\int_0^\infty\phi(\lambda)~dE(\lambda).
\end{equation}
That is, $\phi(L)$ is the operator with domain
$$\Dom(\phi(L))=\set{f\in L^2(\Omega):\int_0^\infty\abs{\phi(\lambda)}^2~dE_{f,f}(\lambda)<\infty},$$
defined by
\begin{equation}\label{def of spectral}
\left\langle\phi(L)f,g\right\rangle_{L^2(\Omega)}=\left\langle\int_0^\infty\phi(\lambda)~dE(\lambda)f,g\right\rangle_{L^2(\Omega)}=\int_0^\infty \phi(\lambda)~dE_{f,g}(\lambda).
\end{equation}
These considerations allow us to define the following operators:
\begin{description}
  \item[The heat-diffusion semigroup generated by $L$] with domain $L^2(\Omega)$,
   $$e^{-tL}=\int_0^\infty e^{-t\lambda}~dE(\lambda),\qquad t\geq0.$$
   We have the contraction property in $L^2(\Omega)$: $\norm{e^{-tL}f}_{L^2(\Omega)}\leq\norm{f}_{L^2(\Omega)}$.
  \item[The fractional operators $L^\sigma$, for $0<\sigma<1$] with domain $\Dom(L^\sigma)\subset\Dom(L)$,
  $$L^\sigma=\int_0^\infty\lambda^\sigma~dE(\lambda)=\frac{1}{\Gamma(-\sigma)}\int_0^\infty\left(e^{-tL}-\iden\right)~\frac{dt}{t^{1+\sigma}}.$$
  \item[The negative powers $L^{-\sigma}$, for $\sigma>0$]
  \begin{equation}\label{eq:L menos sigma}
  L^{-\sigma}=\int_0^\infty\lambda^{-\sigma}~dE(\lambda)=\frac{1}{\Gamma(\sigma)}\int_0^\infty e^{-tL}~\frac{dt}{t^{1-\sigma}}.
  \end{equation}
\end{description}

\begin{proof}[Proof of Theorem \ref{Thm:Extension}]~
\newline\indent\textbf{1.} First we prove that $u(\cdot,y)\in L^2(\Omega)$ and, for all $g\in L^2(\Omega)$,
\begin{equation}\label{u contra g}
\left\langle u(\cdot,y),g(\cdot)\right\rangle_{L^2}=\frac{1}{\Gamma(\sigma)}\int_0^\infty\left\langle e^{-tL}(L^\sigma f),g\right\rangle_{L^2(\Omega)}e^{-\frac{y^2}{4t}}~\frac{dt}{t^{1-\sigma}}.
\end{equation}
For each $R>0$ we let
$$u_R(x,y)=\frac{1}{\Gamma(\sigma)}\int_0^Re^{-tL}(L^\sigma f)(x)e^{-\frac{y^2}{4t}}~\frac{dt}{t^{1-\sigma}}.$$
Since $f\in\Dom(L^\sigma)$, $e^{-tL}(L^\sigma f)\in L^2(\Omega)$. Moreover, $e^{-\frac{y^2}{4t}}/t^{1-\sigma}$ is integrable near $0$ as a function of $t$. Then, using Bochner's Theorem, \eqref{def of spectral}, the fact that $dE_{f,g} (\lambda)$ is of bounded variation, and the change of variables $t=r/\lambda$, we have
\begin{align*}
    \left\langle u_R(\cdot,y),g(\cdot)\right\rangle_{L^2(\Omega)} &= \frac{1}{\Gamma(\sigma)}\int_0^R\left\langle e^{-tL}L^\sigma f,g\right\rangle_{L^2(\Omega)}e^{-\frac{y^2}{4t}}~\frac{dt}{t^{1-\sigma}} \\
    &= \frac{1}{\Gamma(\sigma)}\int_0^R\int_0^\infty e^{-t\lambda}\lambda^\sigma~dE_{f,g}(\lambda)~e^{-\frac{y^2}{4t}}~\frac{dt}{t^{1-\sigma}} \\
    &= \frac{1}{\Gamma(\sigma)}\int_0^\infty\int_0^R e^{-t\lambda}(t\lambda)^\sigma e^{-\frac{y^2}{4t}}~\frac{dt}{t}~dE_{f,g}(\lambda) \\
    &= \frac{1}{\Gamma(\sigma)}\int_0^\infty\int_0^{R\lambda}e^{-r}r^\sigma e^{-\frac{y^2}{4r}\lambda}~\frac{dr}{r}~dE_{f,g}(\lambda),
\end{align*}
so that
$$\abs{\left\langle u_R(\cdot,y),g(\cdot)\right\rangle_{L^2(\Omega)}}\leq\frac{1}{\Gamma(\sigma)}\int_0^\infty\int_0^\infty e^{-r}r^\sigma~\frac{dr}{r}~d\abs{E_{f,g}}(\lambda)\leq\norm{f}_{L^2(\Omega)}\norm{g}_{L^2(\Omega)}.$$
Therefore, for each fixed $y>0$, $u_R(\cdot,y)$ is in $L^2(\Omega)$, and $\norm{u_R(\cdot,y)}_{L^2(\Omega)}\leq\norm{f}_{L^2(\Omega)}$.

The last calculation shows that $\lim_{R_1,R_2\to\infty}\left\langle u_{R_2}(\cdot,y)-u_{R_1}(\cdot,y),g(\cdot)\right\rangle_{L^2(\Omega)}=0$. Then, for any sequence $\set{R^j}_{j\in\mathbb{N}}$ of positive numbers, with $R^j\nearrow\infty$, the family $\set{u_{R^j}(\cdot,y)}_{j\in\mathbb{N}}$ is a Cauchy sequence of bounded linear operators on $L^2(\Omega)$. Thus, $u_R(\cdot,y)\to u(\cdot,y)$ weakly in $L^2(\Omega)$, as $R\to\infty$, and $u(\cdot,y)\in L^2(\Omega)$. Moreover,
\begin{align*}
 \lefteqn{\left\langle u(\cdot,y),g(\cdot)\right\rangle_{L^2(\Omega)}=\lim_{R\to\infty}\left\langle u_R(\cdot,y),g(\cdot)\right\rangle_{L^2(\Omega)}= \lim_{R\to\infty}\frac{1}{\Gamma(\sigma)}\int_0^\infty\int_0^Re^{-t\lambda}(t\lambda)^\sigma e^{-\frac{y^2}{4t}}~\frac{dt}{t}~dE_{f,g}(\lambda)} \\
    &= \frac{1}{\Gamma(\sigma)}\int_0^\infty\int_0^\infty e^{-t\lambda}(t\lambda)^\sigma e^{-\frac{y^2}{4t}}~\frac{dt}{t}~dE_{f,g}(\lambda)=\frac{1}{\Gamma(\sigma)}\int_0^\infty\int_0^\infty e^{-t\lambda}\lambda^\sigma ~dE_{f,g}(\lambda)~e^{-\frac{y^2}{4t}}~\frac{dt}{t^{1-\sigma}} \\
    &= \frac{1}{\Gamma(\sigma)}\int_0^\infty\langle e^{-tL}(L^\sigma f),g\rangle_{L^2(\Omega)}e^{-\frac{y^2}{4t}}~\frac{dt}{t^{1-\sigma}},
\end{align*}
where the limit can be taken inside the integral because the double integral converges absolutely. Hence, \eqref{u contra g} follows.
\newline\indent \textbf{2.} Next we show that  $u(\cdot,y)\in\Dom(L)$, that is,
\begin{equation*}\label{limit}
\lim_{s\to0^+}\left\langle\frac{e^{-sL}u(\cdot,y)-u(\cdot,y)}{s},g(\cdot)\right\rangle_{L^2(\Omega)}\hbox{ exists for all }g\in L^2(\Omega).
\end{equation*}
As $e^{-sL}$ is self-adjoint, by \eqref{u contra g} we have
\begin{align*}
    \left\langle e^{-sL}u(\cdot,y),g(\cdot)\right\rangle_{L^2(\Omega)} &= \left\langle u(\cdot,y),e^{-sL}g(\cdot)\right\rangle_{L^2(\Omega)}=\frac{1}{\Gamma(\sigma)}\int_0^\infty\left\langle e^{-tL}L^\sigma f,e^{-sL}g\right\rangle_{L^2(\Omega)}e^{-\frac{y^2}{4t}}~\frac{dt}{t^{1-\sigma}} \\
    &= \frac{1}{\Gamma(\sigma)}\int_0^\infty\left\langle e^{-sL}e^{-tL}L^\sigma f,g\right\rangle_{L^2(\Omega)}e^{-\frac{y^2}{4t}}~\frac{dt}{t^{1-\sigma}}.
\end{align*}
Hence, \eqref{u contra g}, \eqref{def of spectral}, Fubini's Theorem and dominated convergence give
\begin{align*}
    \left\langle\frac{e^{-sL}u(\cdot,y)-u(\cdot,y)}{s},g(\cdot)\right\rangle_{L^2(\Omega)} &= \frac{1}{\Gamma(\sigma)}\int_0^\infty\left\langle\frac{e^{-sL}e^{-tL}L^\sigma f-e^{-tL}L^\sigma f}{s},g\right\rangle_{L^2(\Omega)}e^{-\frac{y^2}{4t}}~\frac{dt}{t^{1-\sigma}} \\
    &= \frac{1}{\Gamma(\sigma)}\int_0^\infty\int_0^\infty\frac{e^{-s\lambda}e^{-t\lambda}\lambda^\sigma -e^{-t\lambda}\lambda^\sigma}{s}~dE_{f,g}(\lambda)~e^{-\frac{y^2}{4t}}~\frac{dt}{t^{1-\sigma}} \\
    &= \frac{1}{\Gamma(\sigma)}\int_0^\infty\int_0^\infty\frac{e^{-s\lambda}e^{-t\lambda}\lambda^\sigma -e^{-t\lambda}\lambda^\sigma}{s}~e^{-\frac{y^2}{4t}}~\frac{dt}{t^{1-\sigma}}~dE_{f,g}(\lambda) \\
    &\underset{s\to0^+}{\longrightarrow}\frac{1}{\Gamma(\sigma)}\int_0^\infty\int_0^\infty\partial_t (e^{-t\lambda})\lambda^\sigma e^{-\frac{y^2}{4t}}~\frac{dt}{t^{1-\sigma}}~dE_{f,g}(\lambda) \\
    &\qquad =\frac{1}{\Gamma(\sigma)}\int_0^\infty\int_0^\infty\partial_t(e^{-t\lambda}) \lambda^\sigma~dE_{f,g}(\lambda)~e^{-\frac{y^2}{4t}}~\frac{dt}{t^{1-\sigma}} \\
    &\qquad =-\frac{1}{\Gamma(\sigma)}\int_0^\infty\left\langle Le^{-tL}L^\sigma f,g\right\rangle_{L^2(\Omega)}e^{-\frac{y^2}{4t}}~\frac{dt}{t^{1-\sigma}}
\end{align*}
\newline \indent \textbf{3.} We check the boundary condition \eqref{Ext bound value}: for $g\in L^2(\Omega)$, by \eqref{u contra g},
\begin{align*}
    \left\langle u(\cdot,y),g(\cdot)\right\rangle &=\frac1{\Gamma(\sigma)} \int_0^\infty\int_0^\infty e^{-t\lambda}(t\lambda)^\sigma~dE_{f,g}(\lambda)~e^{-\frac{y^2}{4t}}\frac{dt}{t^{1-\sigma}} \\
    &=\frac1{\Gamma(\sigma)} \int_0^\infty\int_0^\infty e^{-r}r^\sigma e^{-\frac{y^2\lambda}{4r}}~\frac{dr}{r}~dE_{f,g}(\lambda)\underset{y\to0}{\longrightarrow}\langle f,g\rangle_{L^2(\Omega)}.
\end{align*}
\newline\indent\textbf{4.} The function $u$ is differentiable with respect to $y$ and
\begin{equation}\label{partial y u}
u_y(x,y)=\frac{1}{\Gamma(\sigma)}\int_0^\infty e^{-tL}(L^\sigma f)(x)~\partial_y(e^{-\frac{y^2}{4t}})~\frac{dt}{t^{1-\sigma}}=\frac{-1}{\Gamma(\sigma)}\int_0^\infty e^{-tL}(L^\sigma f)(x)~\frac{ye^{-\frac{y^2}{4t}}}{2t}~\frac{dt}{t^{1-\sigma}}.
\end{equation}
Indeed applying \eqref{u contra g}, dominated convergence and Bochner's Theorem we get
\begin{align*}
    \lim_{h\to0}\left\langle\frac{u(\cdot,y+h)-u(\cdot)}{h},g(\cdot)\right\rangle_{L^2(\Omega)} &= \frac{1}{\Gamma(\sigma)}\int_0^\infty\left\langle e^{-tL}(L^\sigma f),g\right\rangle_{L^2(\Omega)}\partial_y(e^{-\frac{y^2}{4t}})~\frac{dt}{t^{1-\sigma}} \\
    &= \left\langle\frac{1}{\Gamma(\sigma)}\int_0^\infty e^{-tL}(L^\sigma f)\partial_y(e^{-\frac{y^2}{4t}})~\frac{dt}{t^{1-\sigma}},g\right\rangle_{L^2(\Omega)}.
\end{align*}
\newline\indent\textbf{5.} The function $u$ verifies the extension equation \eqref{Ext equation}. Observe that the integral defining $u_y$ in \eqref{partial y u} is absolutely convergent as a Bochner integral, and it can be differentiated again with respect to $y$. Hence,
\begin{align*}
    \left\langle\frac{1-2\sigma}{y}~u_y(\cdot,y)+u_{yy}(\cdot,y),g(\cdot)\right\rangle_{L^2(\Omega)} &= \frac1{\Gamma(\sigma)}\int_0^\infty\left\langle e^{-tL}L^\sigma f,g\right\rangle_{L^2(\Omega)} \left(\frac{\sigma-1}{t}+\frac{y^2}{4t^2}\right)e^{-\frac{y^2}{4t}}~\frac{dt}{t^{1-\sigma}} \\
    &= -\frac{1}{\Gamma(\sigma)}\int_0^\infty\partial_t\left[\int_0^\infty e^{-t\lambda}\lambda^\sigma~dE_{f,g}(\lambda)\right]e^{-\frac{y^2}{4t}}~\frac{dt}{t^{1-\sigma}} \\
    &= \frac{1}{\Gamma(\sigma)}\int_0^\infty\lambda\int_0^\infty e^{-t\lambda}\lambda^\sigma e^{-\frac{y^2}{4t}}~\frac{dt}{t^{1-\sigma}}~dE_{f,g}(\lambda) \\
    &= \left\langle L\int_0^\infty e^{-tL}(L^\sigma f)e^{-\frac{y^2}{4t}}~\frac{dt}{t^{1-\sigma}},g\right\rangle=\left\langle Lu(\cdot,y),g(\cdot)\right\rangle_{L^2(\Omega)}.
\end{align*}
\newline\indent\textbf{6.} Let us check \eqref{Neumann condtn}. Note that, for all $g\in L^2(\Omega)$, by \eqref{u contra g} and the change of variables $t=y^2/(4r)$,
\begin{align*}
  \left\langle\frac{u(\cdot,y)-u(\cdot,0)}{y^{2\sigma}},g(\cdot)\right\rangle_{L^2(\Omega)}  
  = \frac{1}{4^\sigma\Gamma(\sigma)}\int_0^\infty\left\langle e^{-\frac{y^2}{4r}L}L^\sigma f,g\right\rangle_{L^2(\Omega)}\left(\frac{e^{-r}-1}{r^\sigma}\right)\frac{dr}{r},
\end{align*}
therefore, since $\lim_{t\to0^+}\left\langle e^{-tL}L^\sigma f,g\right\rangle_{L^2(\Omega)}=\left\langle L^\sigma f,g\right\rangle_{L^2(\Omega)}$, by dominated convergence, we obtain the first identity in \eqref{Neumann condtn}. Using \eqref{partial y u} and the same change of variables, the second equality of \eqref{Neumann condtn} follows analogously.
\newline\indent\textbf{7.} Finally, we derive the Poisson formula \eqref{Poisson}. By \eqref{u contra g}, \eqref{def of spectral}, Fubini's Theorem and the change of variables $t=y^2/(4r\lambda)$, we get
\begin{align*}
  \lefteqn{\left\langle u(\cdot,y),g(\cdot)\right\rangle_{L^2(\Omega)}=\frac{1}{\Gamma(\sigma)}\int_0^\infty\int_0^\infty e^{-t\lambda}(t\lambda)^\sigma e^{-\frac{y^2}{4t}}~\frac{dt}{t}~dE_{f,g}(\lambda)} \\
    &\qquad =\frac{1}{\Gamma(\sigma)}\int_0^\infty\int_0^\infty e^{-\frac{y^2}{4r}}\left(\frac{y^2}{4r}\right)^\sigma e^{-r\lambda}~\frac{dr}{r}~dE_{f,g}(\lambda)=\frac{y^{2\sigma}}{4^\sigma\Gamma(\sigma)}\int_0^\infty\left\langle e^{-tL}f,g\right\rangle_{L^2(\Omega)}e^{-\frac{y^2}{4r}}~\frac{dr}{r^{1+\sigma}} \\
    &\qquad =\left\langle \frac{y^{2\sigma}}{4^\sigma\Gamma(\sigma)}\int_0^\infty e^{-tL}f~e^{-\frac{y^2}{4r}}~\frac{dr}{r^{1+\sigma}} ,g\right\rangle_{L^2(\Omega)}.
\end{align*}
The last equality is due to Bochner's Theorem.
\newline\indent The second identity of \eqref{Poisson} follows from the first one via the change of variables $r=y^2/(4t)$.
\end{proof}

In what follows, we assume that the heat-diffusion semigroup is given by integration against a nonnegative \emph{heat kernel} $K_t(x,z)$, $t>0$, $x,z\in\Omega$, that is,
$$e^{-tL}f(x)=\int_\Omega K_t(x,z)f(z)~d\eta(z).$$
Since $e^{-tL}$ is self-adjoint, $K_t(x,z)=K_t(z,x)$. The second assumption we make is that the heat kernel belongs to the domain of $L$, and $\partial_tK_t(x,z)=LK_t(x,z)$, the derivative with respect to $t$ is understood in the classical sense. This implies that
$$\partial_t\int_\Omega K_t(x,z)f(z)~d\eta(z)=\int_\Omega\partial_tK_t(x,z)f(z)~d\eta(z),\qquad f\in L^2(\Omega).$$
Motivated by concrete examples, we add the hypotheses that given $x$, there exists a constant $C_x$ and $\varepsilon>0$, such that $\norm{K_t(x,\cdot)}_{L^2(\Omega)}+\norm{\partial_tK_t(x,\cdot)}_{L^2(\Omega)}\leq C_x(1+t^\varepsilon)t^{-\varepsilon}$.

\begin{thm}[Poisson formula]\label{Thm:Poisson}
Denote by $\mathcal{P}_y^\sigma f(x)$ the function $u(x,y)$ given in \eqref{Poisson}. Then:
\begin{enumerate}
  \item We have $\displaystyle\mathcal{P}_y^\sigma f(x)=\int_\Omega P_y^\sigma(x,z)f(z)~d\eta(z)$, where the \emph{Poisson kernel}
  \begin{equation}\label{Poisson kernel}
    P_y^\sigma(x,z):=\frac{y^{2\sigma}}{4^\sigma\Gamma(\sigma)}\int_0^\infty
    K_t(x,z)e^{-\frac{y^2}{4t}}~\frac{dt}{t^{1+\sigma}},
  \end{equation}
  is, for each fixed $z\in\Omega$, an $L^2(\Omega)$-function that verifies \eqref{Ext equation}.
  \item $\sup_{y\geq0}\abs{\mathcal{P}_y^\sigma f}\leq\sup_{t\geq0}\abs{e^{-tL}f}$, in $\Omega$.
  \item If $e^{-tL}$ has the contraction property in $L^p(\Omega)$, then $\norm{\mathcal{P}_y^\sigma f}_{L^p(\Omega)}\leq\norm{f}_{L^p(\Omega)}$, for all $y\geq0$.
  \item If $\lim_{t\to0^+}e^{-tL}f=f$ in $L^p(\Omega)$, then $\lim_{y\to0^+}\mathcal{P}_y^\sigma f=f$ in $L^p(\Omega)$.
\end{enumerate}
\end{thm}

\begin{proof}
The integral formula in (1)  can be verified by using \eqref{Poisson}, Bochner's and Fubini's Theorems.
\newline \indent In order to see that the Poisson kernel satisfies \eqref{Ext equation}, we begin by showing that it belongs to the domain of $L$. By the assumptions established on the $L^2$-norm of the heat kernel, $P_y^\sigma(\cdot,z)\in L^2(\Omega)$, for each $z$, and, by Bochner's Theorem,
$$e^{-sL}P_y^\sigma(\cdot,z)=\frac{y^{2\sigma}}{4^\sigma\Gamma(\sigma)}\int_0^\infty e^{-sL}K_t(x,z)e^{-\frac{y^2}{4t}}~\frac{dt}{t^{1+\sigma}},\qquad s\geq0.$$
With this, we have
\begin{equation}\label{coc incr}
\frac{e^{-sL}P_y^\sigma(\cdot,z)-P_y^\sigma(\cdot,z)}{s}=\frac{y^{2\sigma}}{4^\sigma\Gamma(\sigma)}\int_0^\infty \frac{e^{-sL}K_t(\cdot,z)-K_t(\cdot,z)}{s}~e^{-\frac{y^2}{4r}}\frac{dr}{r^{1+\sigma}}.
\end{equation}
We use the Mean Value Theorem, the fact that $K_t(\cdot,z)\in\Dom(L)$, and the contraction property of $e^{-sL}$, to get
\begin{align*}
    \norm{\frac{e^{-sL}K_t(\cdot,z)-K_t(\cdot,z)}{s}}_{L^2(\Omega)} &= \norm{Le^{-\theta L}K_t(\cdot,z)}_{L^2(\Omega)}=\norm{e^{-\theta L}LK_t(\cdot,z)}_{L^2(\Omega)} \\
    &\leq \norm{LK_t(\cdot,z)}_{L^2(\Omega)}=\norm{\partial_tK_t(\cdot,z)}_{L^2(\Omega)}\leq C_z(1+t^\varepsilon)t^{-\varepsilon}.
\end{align*}
Hence, the Dominated Convergence Theorem (for Bochner integrals) can be applied in \eqref{coc incr} to see that the limit as $s\to0^+$ of both sides exists, and
$$-L_xP_y^\sigma(x,z)=\frac{y^{2\sigma}}{4^\sigma\Gamma(\sigma)}\int_0^\infty\partial_t\left(K_t(x,z)\right) e^{-\frac{y^2}{4t}}~\frac{dt}{t^{1+\sigma}}.$$
Now we are in position to check that $P_y^\sigma(x,z)$ verifies \eqref{Ext equation}. Note that, by dominated convergence, the derivatives with respect to $y$ of $P_y^\sigma(x,z)$ exist and can be computed by differentiation inside the integral sign in \eqref{Poisson kernel}. Then, using integration by parts,
\begin{align*}
 \frac{1-2\sigma}{y}~\partial_yP_y^\sigma(x,z)+\partial_{yy}P_y^\sigma(x,z) &= \frac{y^{2\sigma}}{4^\sigma\Gamma(\sigma)}\int_0^\infty K_t(x,z)e^{-\frac{y^2}{4t}}\left(\frac{y^2}{4t^2}-\frac{1+\sigma}{t}\right)~\frac{dt}{t^{1+\sigma}} \\
 &= -\frac{y^{2\sigma}}{4^\sigma\Gamma(\sigma)}\int_0^\infty\partial_t(K_t(x,z)) e^{-\frac{y^2}{4t}}~\frac{dt}{t^{1+\sigma}}=L_xP_y^\sigma(x,z),
\end{align*}
thus (1) is proved. (2) follows from the second identity of \eqref{Poisson}. The contraction property of the heat semigroup gives (3):
$$\norm{\mathcal{P}_y^\sigma f}_{L^p(\Omega)}\leq\frac{1}{\Gamma(\sigma)}\int_0^\infty \big\|e^{-\frac{y^2}{4r}L}f\big\|_{L^p(\Omega)}e^{-r}~\frac{dr}{r^{1-\sigma}}\leq\norm{f}_{L^p(\Omega)}.$$
Finally, observe that
$$\norm{\mathcal{P}_y^\sigma f-f}_{L^p(\Omega)}\leq\frac{1}{\Gamma(\sigma)}\int_0^\infty\big\|e^{-\frac{y^2}{4t}L}f-f\big\|_{L^p(\Omega)}e^{-r}~\frac{dr}{r^{1-\sigma}},$$
so (4) follows.
\end{proof}

\begin{rem}
Note in \eqref{Poisson} that, when $\sigma=1/2$, $\mathcal{P}_y^{1/2}f=e^{-t\sqrt{L}}f$ is the classical subordinated Poisson semigroup of $L$ acting on $f$ (see \cite[p.~47~and~49]{SteinTopics}).
\end{rem}

\begin{prop}[Fundamental solution of \eqref{Ext equation}]\label{Prop:Fund sol}
The function
\begin{equation}\label{eq:Fund sol}
\Psi_x^\sigma(z,y)=\frac{1}{\Gamma(\sigma)}\int_0^\infty K_t(x,z)e^{-\frac{y^2}{4t}}~\frac{dt}{t^{1-\sigma}},
\end{equation}
satisfies equation \eqref{Ext equation}, $\Psi_x^\sigma(z,y)=\Psi_z^\sigma(x,y)$, and
\begin{equation}\label{Fund sol to0}
\lim_{y\to0^+}\left\langle\frac{1}{2\sigma}~y^{1-2\sigma}\partial_y\Psi_x^\sigma(\cdot,y),f(\cdot)\right\rangle_{L^2(\Omega)}= \frac{\Gamma(-\sigma)}{4^\sigma\Gamma(\sigma)}f(x).
\end{equation}
\end{prop}

\begin{proof}
As in the proof of Theorem \ref{Thm:Extension}, it can be checked that for each $x$,
\begin{align*}
    \lim_{R\to\infty}\left\langle\frac{1}{\Gamma(\sigma)}\int_0^R K_t(x,\cdot)e^{-\frac{y^2}{4t}}~\frac{dt}{t^{1-\sigma}},g(\cdot)\right\rangle_{L^2(\Omega)} &= \left\langle\Psi_x^\sigma(\cdot,y),g(\cdot)\right\rangle_{L^2(\Omega)} \\
    &= \frac{1}{\Gamma(\sigma)}\int_0^\infty e^{-tL}g(x)e^{-\frac{y^2}{4t}}~\frac{dt}{t^{1-\sigma}},
\end{align*}
and that \eqref{eq:Fund sol} satisfies \eqref{Ext equation}. Differentiation with respect to $y$ inside the integral in \eqref{eq:Fund sol} can be performed to get
$$\frac{y^{1-2\sigma}}{2\sigma}~\partial_y\Psi_x^\sigma(z,y)=\frac{-1}{4^\sigma\sigma\Gamma(\sigma)}\int_0^\infty K_t(x,z)e^{-\frac{y^2}{4t}}\left(\frac{y^2}{4t}\right)^{1-\sigma}\frac{dt}{t}=\frac{-1}{4^\sigma\sigma\Gamma(\sigma)}\int_0^\infty K_{\frac{y^2}{4r}}(x,z)e^{-r}~\frac{dr}{r^\sigma}.$$
With this we obtain \eqref{Fund sol to0}:
\begin{align*}
  \frac{y^{1-2\sigma}}{2\sigma}\int_\Omega\partial_y\Psi_x^\sigma(z,y)f(z)~d\eta(z) &= \frac{-1}{4^\sigma\sigma\Gamma(\sigma)}\int_0^\infty e^{-\frac{y^2}{4r}}f(x)e^{-r}~\frac{dr}{r^\sigma} \\
   & \to\frac{-\Gamma(1-\sigma)}{4^\sigma\sigma\Gamma(\sigma)}f(x)=\frac{\Gamma(-\sigma)}{4^\sigma\Gamma(\sigma)}f(x),\qquad y\to0^+.
\end{align*}
\end{proof}

\begin{rem}
It can also be proved that
$$\lim_{y\to0^+}\left\langle\frac{\Psi_x^\sigma(\cdot,y)-\Psi_x^\sigma(\cdot,0)}{y^{2\sigma}},f(\cdot)\right\rangle_{L^2(\Omega)} =\frac{\Gamma(-\sigma)}{4^\sigma\Gamma(\sigma)}f(x).$$
\end{rem}

\begin{prop}
Let $v(x,y)=y^{1-2\sigma}u_y(x,y)$, where $u$ solves \eqref{Ext equation}. Then $v$ is a solution of the following ``conjugate equation''
\begin{equation}\label{Conj L}
-Lv-\frac{1-2\sigma}{y}~v_y+v_{yy}=0,\qquad\hbox{ in }\Omega\times(0,\infty).
\end{equation}
\end{prop}

\begin{proof}
The calculation is analogous to the one given in \cite{Caffarelli-Silvestre CPDE}, with the obvious modifications, and we omit it here.
\end{proof}

\begin{rem}\label{Rem:Sol fund y Poisson}
As in \cite{Caffarelli-Silvestre CPDE} the fundamental solution \eqref{eq:Fund sol} and the ``conjugate equation'' \eqref{Conj L} (which coincides with the conjugate equation given in \cite{Caffarelli-Silvestre CPDE} when $L=-\Delta$) can help us to find the Poisson kernel \eqref{Poisson kernel}. Indeed, we want to write $u(x,y)=\mathcal{P}_y^\sigma f(x)=\int_\Omega P_y^\sigma(x,z)f(z)~d\eta(z)$ where the Poisson kernel $P_y^\sigma(x,z)$ must be a solution of \eqref{Ext equation} for all $z$ and $\lim_{y\to0^+}\mathcal{P}_y^\sigma f(x)=f(x)$. The right choice would be
\begin{equation}\label{Poisson y Fund sol}
P_y^\sigma(x,z)=\frac{4^{1-\sigma}\Gamma(1-\sigma)}{\Gamma(-(1-\sigma))2(1-\sigma)}~y^{1-2(1-\sigma)}\partial_y\Psi_x^{1-\sigma}(z,y)=C_{1-\sigma}y^{1-2(1-\sigma)}\partial_y\Psi_x^{1-\sigma}(z,y),
\end{equation}
since it solves the ``conjugate equation'' \eqref{Conj L} with $1-\sigma$ in the place of $\sigma$ (thus it verifies \eqref{Ext equation}) and by \eqref{Fund sol to0} and the choice of $C_{1-\sigma}$,
$$\lim_{y\to0^+}C_{1-\sigma}\int_\Omega y^{1-2(1-\sigma)}\partial_y\Psi_x^{1-\sigma}(z,y)f(z)~d\eta(z)=f(x).$$
A simple calculation shows that \eqref{Poisson y Fund sol} coincides with \eqref{Poisson kernel}.
\end{rem}

In the following discussion we shall assume that the operator $L$ can be factorized as
$L=D_i^\ast D_i$, where $D_i=a_i(x_i)\partial_{x_i}+b_i(x_i)$, is a one dimensional (in the $i$th direction) partial differential operator, and $D_i^\ast$ is the formal adjoint (with respect to $d\eta$) of $D_i$. See examples at the end of this Section. In this case we give a definition of $n$ conjugate functions related to the Poisson formula for $u$.

Let $E_\sigma=-L+\frac{1-2\sigma}{y}\partial_y+\partial_{yy}$. Then the factorization
$$E_\sigma=-\sum_{i=1}^nD_i^\ast D_i+y^{-(1-2\sigma)}\partial_y(y^{1-2\sigma}\partial_y),$$
suggests the following definition of Cauchy-Riemann equations for a system of functions $u,v_1,\ldots,v_n:\Omega\times(0,\infty)\to\Real$ such that $E_\sigma u=0$:
\begin{equation}\label{CR Extension}
\left\{
\begin{aligned}
y^{1-2\sigma}\partial_yu&=D_1^\ast v_1+\cdots+D_n^\ast v_n,\\
D_iu&=y^{-(1-2\sigma)}\partial_yv_i,&&i=1,\ldots,n,\\
D_kv_i&=D_iv_k,&&i,k=1,\ldots,n.
\end{aligned}
\right.
\end{equation}

\begin{prop}
Let $u$ be a solution of $E_\sigma u=0$ in $\Omega\times(0,\infty)$. If $v_1,\ldots,v_n$ verify \eqref{CR Extension} then each $v_i$ solves the \emph{$i$th conjugate equation}
\begin{equation}\label{Conj Li}
E^i_{1-\sigma}v_i=-Lv_i+[D_i^\ast,D_i]v_i-\frac{1-2\sigma}{y}~\partial_yv_i+\partial_{yy}v_i=0,\qquad i=1,\ldots,n,
\end{equation}
where $[D_i^\ast,D_i]=D_i^\ast D_i-D_iD_i^\ast$.
\end{prop}

\begin{proof}
\begin{align*}
  -Lv_i+[D_i^\ast,D_i]v_i &= -\sum_{k\neq i}D_k^\ast D_kv_i-D_iD_i^\ast v_i=-\sum_{k\neq i}D_k^\ast D_iv_k-D_iD_i^\ast v_i=-D_i\left(\sum_{k=1}^nD_k^\ast v_k\right) \\
   &= -D_i\left(y^{1-2\sigma}\partial_yu\right)=-y^{1-2\sigma}\partial_y(D_iu)=-y^{1-2\sigma}\partial_y(y^{-(1-2\sigma)}\partial_yv_i) \\
   &= \frac{1-2\sigma}{y}~\partial_yv_i-\partial_{yy}v_i.
\end{align*}
\end{proof}

\begin{rem}
The $i$th conjugate equation \eqref{Conj Li} is not the same as the ``conjugate equation'' \eqref{Conj L}. They will coincide only when $[D_i^\ast,D_i]=0$. This is the case if $L=-\Delta$: the conjugate equation  established in \cite{Caffarelli-Silvestre CPDE} is equal to each $i$th conjugate equation \eqref{Conj Li}.
\end{rem}

\begin{prop}\label{Prop:Conj Poisson kernel}
Fix $z\in\Omega$ and choose $u(x,y)=P_y^\sigma(x,z)$. Then a solution to \eqref{CR Extension} is given by the $n$ \emph{conjugate Poisson kernels} defined by
\begin{equation}\label{Poisson kernel conj}
v_i(x,y):=Q_y^{\sigma,i}(x,z)=\frac{-2}{4^\sigma\Gamma(\sigma)}D_i\int_0^\infty K_t(x,z)e^{-\frac{y^2}{4t}}~\frac{dt}{t^\sigma},\qquad i=1,\ldots,n.
\end{equation}
\end{prop}

\begin{proof}
From \eqref{Poisson y Fund sol} and the second equation of \eqref{CR Extension} we have $C_{1-\sigma}D_i\partial_y\Psi_x^{1-\sigma}(z,y)=\partial_yQ_y^{\sigma,i}(x,z)$, so, in view of \eqref{eq:Fund sol}, $Q_y^{\sigma,i}(x,z)$ can be chosen as in \eqref{Poisson kernel conj}. Clearly $D_kQ_y^{\sigma,i}(x,z)=D_iQ_y^{\sigma,k}(x,z)$. Moreover, the first equation of \eqref{CR Extension} also holds:
\begin{align*}
  D_1^\ast Q_y^{\sigma,1}(x,z)+\cdots+D_n^\ast Q_y^{\sigma,n}(x,z) &= \frac{-2}{4^\sigma\Gamma(\sigma)}\sum_{i=1}^nD_i^\ast D_i\int_0^\infty K_t(x,z)e^{-\frac{y^2}{4t}}~\frac{dt}{t^\sigma} \\
   &= \frac{-2}{4^\sigma\Gamma(\sigma)}\int_0^\infty LK_t(x,z)e^{-\frac{y^2}{4t}}~\frac{dt}{t^\sigma} \\
   &= \frac{2}{4^\sigma\Gamma(\sigma)}\int_0^\infty \partial_tK_t(x,z)e^{-\frac{y^2}{4t}}~\frac{dt}{t^\sigma} \\
   &= \frac{-2}{4^\sigma\Gamma(\sigma)}\int_0^\infty K_t(x,z)e^{-\frac{y^2}{4t}}\left(\frac{y^2}{4t}-\sigma\right)~\frac{dt}{t^{1+\sigma}} \\
   &= y^{1-2\sigma}\partial_yP_y^\sigma(x,z).
\end{align*}
\end{proof}

\begin{cor}
The \emph{Poisson integral of $f$}, $u(x,y)=\mathcal{P}_y^\sigma f(x)$, and the $n$ \emph{conjugate Poisson integrals} of $f$ defined by
\begin{equation}\label{Poisson conjugado}
v_i(x,y)\equiv\mathcal{Q}_y^{\sigma,i}f(x):=\int_\Omega Q_y^{\sigma,i}(x,z)f(z)~d\eta(z)=\frac{-2}{4^\sigma\Gamma(\sigma)}D_i\int_0^\infty e^{-tL}f(x)e^{-\frac{y^2}{4t}}~\frac{dt}{t^\sigma},
\end{equation}
for $i=1,\ldots,n$, solve \eqref{CR Extension}.
\end{cor}

\begin{rem}
When $\sigma=1/2$, $\mathcal{Q}_y^{1/2,i}f(x)$ is the $i$th conjugate function of $f$ associated to $L$, see \cite{SteinTopics} and \cite{Thangavelu}. A natural question arises: what is the limit of $\mathcal{Q}_y^{\sigma,i}f(x)$ as $y\to0^+$? The answer is contained in the next result.
\end{rem}

\begin{thm}\label{Thm:Conj Poisson to 0}
For each $x\in\Omega$,
$$\lim_{y\to0^+}\mathcal{Q}_y^{\sigma,i}f(x)=\frac{-2\Gamma(1-\sigma)}{4^\sigma\Gamma(\sigma)}D_iL^{-(1-\sigma)}f(x).$$
\end{thm}

\begin{proof}
From the expression of $\mathcal{Q}_y^{\sigma,i}f(x)$ in \eqref{Poisson conjugado} and \eqref{eq:L menos sigma},
$$\lim_{y\to0^+}\mathcal{Q}_y^{\sigma,i}f(x)=\frac{-2\Gamma(1-\sigma)}{4^\sigma\Gamma(\sigma)}D_i\frac{1}{\Gamma(1-\sigma)}\int_0^\infty e^{-tL}f(x)\frac{dt}{t^{1-(1-\sigma)}}=\frac{-2\Gamma(1-\sigma)}{4^\sigma\Gamma(\sigma)}D_iL^{-(1-\sigma)}f(x).$$
\end{proof}

\begin{rem}
The conclusion of Theorem \ref{Thm:Conj Poisson to 0} can also be obtained from the following observation: except for a multiplicative constant, the last formula of \eqref{Poisson conjugado} is just the $D_i$-derivative of the solution of the extension problem \eqref{Ext equation} for $L^{1-\sigma}$ with boundary value $L^{-(1-\sigma)}f(x)$ (see \eqref{Extension}). For $\sigma=1/2$, Theorem \ref{Thm:Conj Poisson to 0} establishes the boundary convergence to the Riesz transforms $D_iL^{-1/2}$ (which in case $L=-\Delta$ are the classical Riesz transforms $\partial_{x_i}(-\Delta)^{-1/2}$). See \cite{SteinTopics} and \cite{Thangavelu}.
\end{rem}

\begin{exams}\label{Examples}
We present now some examples of operators $L$ for which our results apply.
\begin{description}
  \item[The Laplacian in $\Real^n$] Observe that, when $f\in\mathcal{S}$,
  \begin{equation}\label{eq:heat Lap}
    e^{t\Delta}f(x)=\int_{\Real^n}W_t(x-z)f(z)~dz,\quad W_t(x)=\frac{1}{(4\pi t)^{n/2}}~e^{-\frac{\abs{x}^2}{4t}}.
  \end{equation}
 The Poisson formula given in \cite{Caffarelli-Silvestre CPDE} is recovered: use the change of variables $\frac{\abs{x-z}^2+y^2}{4t}=r$, in \eqref{Poisson kernel},Ê to see that the Poisson kernel in this case is $$P_y^{\sigma,-\Delta}(x,z)=\frac{y^{2\sigma}}{4^\sigma\Gamma(\sigma)}\int_0^\infty\frac{e^{-\frac{\abs{x-z}^2 +y^2}{4t}}}{(4\pi t)^{n/2}}~\frac{dt}{t^{1+\sigma}}=\frac{\Gamma(n/2+\sigma)}{\pi^{n/2}\Gamma(\sigma)}\cdot \frac{y^{2\sigma}}{\left(\abs{x-z}^2+y^2\right)^{\frac{n+2\sigma}{2}}}.$$
  The function $P_y^{1/2,-\Delta}(x,z)$ is the classical Poisson kernel for the harmonic extension of a function to the upper half space. The Cauchy-Riemann equations read
    \begin{equation}\label{CR Laplacian}
    \left\{
        \begin{aligned}
            y^{1-2\sigma}\partial_yu&=-\left(\partial_{x_1}v_1+\cdots+\partial_{x_n}v_n\right),\\
            \partial_{x_i}u&=y^{-(1-2\sigma)}\partial_yv_i,&&i=1,\ldots,n,\\
            \partial_{x_k}v_i&=\partial_{x_i}v_k,&&i,k=1,\ldots,n.
        \end{aligned}
    \right.
    \end{equation}
  The case $\sigma=1/2$ is the classical Cauchy-Riemann system for the $n$ conjugate harmonic functions to $u$. In dimension one \eqref{CR Laplacian} reduces to
    \begin{equation*}
    \left\{
        \begin{aligned}
            y^{1-2\sigma}\partial_yu&=-\partial_x v,\\
            \partial_xu&=y^{-(1-2\sigma)}\partial_yv,
        \end{aligned}
    \right.
    \end{equation*}
    which already appeared in \cite{Muckenhoupt-Stein}.
      \item[Classical expansions] $L$ can be each one of the operators arising in orthogonal expansions, like the Ornstein-Uhlenbeck operator (Hermite polynomials and Gaussian measure $d\eta(x)=e^{-\abs{x}^2}~dx$),
      $$-\Delta+2x\cdot\nabla=\sum_i\left(-\partial_{x_i}+2x_i\right)\left(\partial_{x_i}\right);$$
      the harmonic oscillator (Hermite functions and Lebesgue measure $d\eta(x)=dx$),
      $$-\Delta+\abs{x}^2 =\frac{1}{2}\sum_i\left[\left(-\partial_{x_i}+x_i\right)\left(\partial_{x_i}+x_i\right)+ \left(\partial_{x_i}+x_i\right)\left(-\partial_{x_i}+x_i\right)\right];$$
      the Laguerre operator (Laguerre polynomials and measure $d\eta(x)=\prod_ix_i^{\alpha_i}e^{-x_i}$)
      $$\sum_ix_i\partial^2_{x_i,x_i}+(\alpha_i+1-x_i)\partial_{x_i}=\sum_i\sqrt{x_i}\left(\partial_{x_i}+ \left(\frac{\alpha_i+1/2}{x_i}-1\right)\right)\sqrt{x_i}~\partial_{x_i};$$
      Jacobi and ultraspherical on $(-1,1)$; etc.
      \newline \indent We would like to point out that in these cases, due to the existence of  smooth eigenfunctions, the proof of Theorem \ref{Thm:Extension} can be performed as an exercise of convergence of orthogonal systems, and it makes it technically simpler.
  \item[Elliptic operators]  Let $L$ be a positive self-adjoint linear elliptic partial differential operator on $L^2(\Omega)$, with Dirichlet boundary conditions, and bounded measurable coefficients. Then the heat kernel exists, and it verifies our assumptions stated before Theorem  \ref{Thm:Poisson}. Even more, its heat kernel has Gaussian bounds \cite[p.~89]{Davies}. We can also consider Schr\"odinger operators with nonnegative potentials in a large class \cite[Section~4.5]{Davies}.
\end{description}
\end{exams}

\section{Existence and uniqueness results for the extension problem}\label{Section Unicidad}

In this section we derive the concrete solution of the extension problem in the case of discrete  spectrum. We also find solutions with null Neumann condition. This is done by using classical Fourier's method.
\newline\indent Let $\set{\phi_k}_{k\in\mathbb{N}_0}$ be an orthonormal basis of $L^2(\Omega)$ such that $L\phi_k=\lambda_k\phi_k$, $\lambda_k\geq0$. Recall the definition of $L^\sigma $ given in \eqref{en L2}.

\subsection{$L^2$ theory}

Let $f\in L^2(\Omega)$ and look for solutions $u$ to \eqref{Ext bound value}-\eqref{Ext equation} of the form
\begin{equation}\label{serie}
u(x,y)=\sum_kc_k(y)\phi_k(x).
\end{equation}
Then for each $k\geq0$ we have to solve the following ordinary differential equation:
$$-\lambda_kc_k+\frac{1-2\sigma}{y}~c_k'+c_k''=0,\qquad y>0,$$
with initial condition $c_k(0)=\langle f,\phi_k\rangle$. According to \cite[p.~106]{Lebedev}, this last equation has a general solution of the form
\begin{equation}\label{general}
c_k(y)=y^\sigma Z_\sigma(\pm i\lambda_k^{1/2}y),
\end{equation}
where $Z_\sigma$ is a linear combination of Bessel functions of order $\sigma$. To have uniqueness of the solution include the boundary condition $\lim_{y\to\infty}u(x,y)=0$, weakly in $L^2(\Omega)$, which translates to the coefficients as
\begin{equation}\label{condicion}
\lim_{y\to\infty}c_k(y) =0.
\end{equation}
From \cite[p.~104]{Lebedev}, $Z_\sigma$ can be written as
\begin{eqnarray}
  \label{AB} Z_\sigma(z) &=& A_1J_\sigma(z)+A_2H_\sigma^{(1)}(z)~=~B_1J_\sigma(z)+B_2H_\sigma^{(2)}(z) \\
  \label{CD} &=& C_1J_\sigma(z)+C_2J_{-\sigma}(z)~=~D_1H_\sigma^{(1)}(z)+D_2H_\sigma^{(2)}(z),
\end{eqnarray}
where $J_\sigma$ denotes the Bessel function of the first kind and $H_\sigma^{(1)}$ and $H_\sigma^{(2)}$ are the Hankel functions. To fulfill condition \eqref{condicion} we need to review the asymptotic behavior of the Bessel functions. When $\abs{\arg z}\leq\pi-\delta$,
\begin{eqnarray}
  \label{J infty} J_\sigma(z) &=& \left(\frac{2}{\pi z}\right)^{1/2}\left[\cos\left(z-\frac{2\sigma\pi+\pi}{4}\right)\left(1+O(\abs{z}^{-2})\right)\right. \\
  \nonumber &&\qquad\qquad\quad\left.-\sin\left(z-\frac{2\sigma\pi+\pi}{4}\right)\left(\frac{4\sigma^2-1}{8z}+O(\abs{z}^{-3})\right)\right], \\
  \nonumber H_\sigma^{(1)}(z) &=& \left(\frac{2}{\pi z}\right)^{1/2}e^{i\left(z-\frac{2\sigma\pi+\pi}{4}\right)}\left(1+O(\abs{z}^{-1})\right), \\
  \nonumber H_\sigma^{(2)}(z) &=& \left(\frac{2}{\pi z}\right)^{1/2}e^{-i\left(z-\frac{2\sigma\pi+\pi}{4}\right)}\left(1+O(\abs{z}^{-1})\right).
\end{eqnarray}
Note that for purely imaginary $z$, $J_\sigma(z)\to\infty$ exponentially and $H_\sigma^{(1)}(z)\to0$ or $\infty$ depending on the sign of the imaginary part of $z$. Putting $z=i\lambda_k^{1/2}y$ in \eqref{general} we see that the only possible choice as solution is the first linear combination of \eqref{AB} as soon as $A_1=0$: $c_k(y)=A_{2,k}y^\sigma H_\sigma^{(1)}(i\lambda_k^{1/2}y)$. If $K_\sigma$ denotes the modified Bessel function of the third kind then $H_\sigma^{(1)}(iz)=2\pi^{-1}i^{-\sigma-1}K_\sigma(z)$ and $$c_k(y)=A_{2,k}y^\sigma\frac{2i^{-\sigma-1}}{\pi}~K_\sigma(\lambda_k^{1/2}y).$$
To determine $A_{2,k}$ use the initial condition. The asymptotic behavior of $K_\sigma(z)$ as $z\to0$ reads
\begin{equation}\label{k en cero}
K_\sigma(z)\approx\Gamma(\sigma)2^{\sigma-1}\frac{1}{z^\sigma}.
\end{equation}
So that, when $y\to0$, $c_k(y)\approx A_{2,k}2^\sigma\pi^{-1}i^{-\sigma-1}\Gamma(\sigma)\lambda_k^{-\sigma/2}$. Therefore
$$A_{2,k}=\frac{\pi i^{1+\sigma}}{2^\sigma\Gamma(\sigma)}\lambda_k^{\sigma/2}\langle f,\phi_k\rangle.$$
Thus
\begin{equation}\label{ck}
c_k(y)=y^\sigma\frac{2^{1-\sigma}}{\Gamma(\sigma)}\lambda_k^{\sigma/2}\langle f,\phi_k\rangle K_\sigma(\lambda_k^{1/2}y).
\end{equation}
Since as $\abs{z}\to\infty$,
$$K_\sigma(z)\approx\left(\frac{\pi}{2z}\right)^{1/2}e^{-z}\left(1+O(\abs{z}^{-1})\right),$$
the series in \eqref{serie}, with $c_k$ as in \eqref{ck}, converges in $L^2(\Omega)$ for each $y\in(0,\infty)$.
 Finally, \eqref{k en cero} implies that \eqref{Ext bound value} is fulfilled in the $L^2(\Omega)$ sense.

On the other hand, by using the properties of the derivatives of $K_\sigma$ (see \cite[p.~110]{Lebedev}) and \eqref{k en cero}, as $y\to0$ we have
\begin{align*}
  \frac{1}{2\sigma}y^{1-2\sigma}c_k'(y) &= \frac{1}{2\sigma}y^{1-2\sigma}\frac{2^{1-\sigma}}{\Gamma(\sigma)}\langle f,\phi_k\rangle\frac{d}{d(\lambda_k^{1/2}y)}\left[(\lambda_k^{1/2}y)^\sigma K_\sigma(\lambda_k^{1/2}y)\right]\frac{d(\lambda_k^{1/2}y)}{y} \\
   &= \frac{1}{2\sigma}y^{1-2\sigma}\frac{2^{1-\sigma}}{\Gamma(\sigma)}\langle f,\phi_k\rangle(-1)(\lambda_k^{1/2}y)^\sigma K_{\sigma-1}(\lambda_k^{1/2}y)\lambda_k^{1/2} \\
   &= \frac{2^{-\sigma}}{-\sigma\Gamma(\sigma)}\langle f,\phi_k\rangle\lambda_k^{\sigma/2}\lambda_k^{1/2}y^{1-\sigma}K_{1-\sigma}(\lambda_k^{1/2}y) \\
   &\approx \frac{2^{-\sigma}}{-\sigma\Gamma(\sigma)}\langle f,\phi_k\rangle\lambda_k^{\sigma/2}\lambda_k^{1/2}y^{1-\sigma}\Gamma(1-\sigma)2^{-\sigma}\frac{1}{(\lambda_k^{1/2}y)^{1-\sigma}} \\
   &= \frac{2^{-2\sigma}\Gamma(1-\sigma)}{-\sigma\Gamma(\sigma)}\lambda_k^\sigma\langle f,\phi_k\rangle=\frac{\Gamma(-\sigma)}{4^\sigma\Gamma(\sigma)}\lambda_k^\sigma\langle f,\phi_k\rangle.
\end{align*}
As a consequence,
$$\frac{1}{2\sigma}\lim_{y\to0^+}y^{1-2\sigma}u_y(x,y)=\frac{\Gamma(-\sigma)}{4^\sigma\Gamma(\sigma)}\sum_k\lambda_k^\sigma\langle f,\phi_k\rangle\phi_k(x)=\frac{\Gamma(-\sigma)}{4^\sigma\Gamma(\sigma)}L^\sigma f(x),$$
the limit taken in the $L^2(\Omega)$-sense (see \eqref{en L2}).


\subsection{Local Neumann solutions}\label{Subsection Neumann}

Let us find a solution to \eqref{Ext equation} such that
\begin{equation}\label{condicion2}
\frac{1}{2\sigma}\lim_{y\to0^+}y^{1-2\sigma}u_y(x,y)=0,\qquad\hbox{ for all }x\in\Omega.
\end{equation}
Writing $u(x,y)=\sum_kd_k(y)\phi_k(x)$, condition \eqref{condicion2} implies that $\lim_{y\to0^+}y^{1-2\sigma}d_k'(y)=0$. Therefore as (see \eqref{general})
$$d_k'(y)=(i\lambda_k)^{1-\sigma}\frac{d}{d(i\lambda_k^{1/2}y)}\left[(i\lambda_k^{1/2}y)^\sigma Z_\sigma(i\lambda_k^{1/2}y)\right]=i\lambda_k^{1/2}y^\sigma Z_{\sigma-1}(i\lambda_k^{1/2}y),$$
we require
$$y^{1-2\sigma}d_k'(y)=i\lambda_k^{1/2}y^{1-\sigma}Z_{\sigma-1}(i\lambda_k^{1/2}y)\to0,\qquad y\to0.$$
When $z\to0$ (see \cite{Lebedev}),
$$J_\sigma(z)\approx \frac{z^\sigma}{2^\sigma\Gamma(1+\sigma)},\quad H_\sigma^{(1)}(z)\approx\frac{2^\sigma\Gamma(\sigma)}{i\pi}\frac{1}{z^\sigma},\quad\hbox{and}\quad H_\sigma^{(2)}(z)\approx-\frac{2^\sigma\Gamma(\sigma)}{i\pi}\frac{1}{z^\sigma}.$$
Then, as $y\to0$,
\begin{align*}
  y^{1-\sigma}J_{\sigma-1}(i\lambda_k^{1/2}y) &\to \frac{(i\lambda_k^{1/2})^{\sigma-1}}{2^{\sigma-1}\Gamma(\sigma)}, \\
  y^{1-\sigma}H_{\sigma-1}^{(1)}(i\lambda_k^{1/2}y)=y^{1-\sigma}i^{2\sigma}H_{1-\sigma}^{(1)}(i\lambda_k^{1/2}y) &\to \frac{2^{1-\sigma}\Gamma(1-\sigma)i^{2\sigma-1}}{\pi}~(i\lambda_k^{1/2})^{1-\sigma}, \\
  y^{1-\sigma}H_{\sigma-1}^{(2)}(i\lambda_k^{1/2}y)=y^{1-\sigma}i^{-2\sigma}H_{1-\sigma}^{(2)}(i\lambda_k^{1/2}y) &\to -\frac{2^{1-\sigma}\Gamma(1-\sigma)i^{-(2\sigma+1)}}{\pi}~(i\lambda_k^{1/2})^{1-\sigma},
\end{align*}
but
$$y^{1-\sigma}J_{1-\sigma}(i\lambda_k^{1/2}y)\approx\frac{(i\lambda_k^{1/2})^{1-\sigma}}{2^{1-\sigma}\Gamma(2-\sigma)}~y^{2-2\sigma}\to0.$$
Consequently, choose the first linear combination in \eqref{CD} with $C_1=0$. Thus
$$d_k(y)=C_{2,k}y^\sigma J_{-\sigma}(i\lambda_k^{1/2}y),$$
verifies $\lim_{y\to0}y^{1-2\sigma}d_k'(y)=0$. So $u$ formally reads
$$u(x,y)=y^\sigma\sum_kC_{2,k}J_{-\sigma}(i\lambda_k^{1/2}y)\phi_k(x).$$
In order to have a convergent series (at least for small $y$), let us determine $C_{2,k}$. Taking into account \eqref{J infty} it is enough to fix $R>0$ and put $C_{2,k}=Ce^{-\lambda_k^{1/2}R}$.

In this way we obtained a solution $u$ to equation \eqref{Ext equation} in $\Omega\times(0,R)$ that satisfies the required property \eqref{condicion2}.

\section{The Harnack's inequality for $H^\sigma$}\label{Section Harnack}

To prove the Harnack's inequality in Theorem \ref{Thm:Harnack} we first study the problem \eqref{Ext bound value}-\eqref{Ext equation} for the harmonic oscillator $L=H=-\Delta+\abs{x}^2$ posed in $\Real^n$ with the Lebesgue measure $d\eta=dx$.

Fix $1\leq p<\infty$ and $N>0$. Define the space
\begin{equation}\label{LpN}
L^p_N=\set{u:\Real^n\to\Real:\norm{u}_{L^p_N}=\left(\int_{\Real^n}\frac{\abs{u(z)}^p}{(1+\abs{z}^2)^{Np}}~dz\right)^{1/p}<\infty}.
\end{equation}
Then $L^p_N\subset\mathcal{S}'$.

The heat semigroup generated by $H$ (see \cite{Thangavelu}) can be given as an integral operator
\begin{equation}\label{eq:heat H}
e^{-tH}f(x)=\int_{\Real^n}G_t(x,z)f(z)~dz=\int_{\Real^n}\frac{e^{-\left[\frac{1}{2}\abs{x-z}^2\coth2t+x\cdot z\tanh t\right]}}{(2\pi\sinh2t)^{n/2}}~f(z)~dz.
\end{equation}
We collect some useful facts about $e^{-tH}$ in the next Proposition, whose proof is postponed to the end of Section \ref{Section Hermite}.

\begin{prop}\label{Prop:Heat LpN}
For $f\in L^p_N$, the heat semigroup $e^{-tH}f(x)$ is well defined and
\begin{equation}\label{calor LpN}
\abs{e^{-tH}f(x)}\leq C\frac{(1+\abs{x}^\rho)\norm{f}_{L^p_N}}{t^{n/2}},\qquad x\in\Real^n,~t>0,
\end{equation}
where $\rho>0$ depends on $p$ and $N$. Moreover, $(\partial_t+H)e^{-tH}f(x)=0$ for all $x\in\Real^n$ and $t>0$ and for $i,j=1,\ldots,n$,
\begin{equation}\label{dentro}
\abs{\partial_{x_i}(e^{-tH}f)(x)}\leq C\frac{(1+\abs{x}^\rho) \norm{f}_{L^p_N}}{t^{(n+1)/2}},\qquad\abs{\partial_{x_ix_j}(e^{-tH}f)(x)}\leq C\frac{(1+\abs{x}^\rho) \norm{f}_{L^p_N}}{t^{(n+2)/2}}.
\end{equation}
If $f$ is also a $C^2$ function in some open subset ${\mathcal{O}}\subset\Real^n$ then $\lim_{t\to0}e^{-tH}f(x)=f(x)$ for all $x\in{\mathcal{O}}$.
\end{prop}

In the particular case we are considering in this Section, Theorem \ref{Thm:Extension} takes the following form, in which the relevant observation is that all identities are classical.

\begin{thm}\label{Thm:Extension smooth}
If $f\in L^p_N$ is a $C^2$ function in some open subset ${\mathcal{O}}\subset\Real^n$ then
\begin{equation}\label{Poisson H y Delta}
u(x,y):=\frac{y^{2\sigma}}{4^\sigma\Gamma(\sigma)}\int_0^\infty e^{-tH}f(x)e^{-\frac{y^2}{4t}}~\frac{dt}{t^{1+\sigma}},
\end{equation}
is well defined for all $x\in\Real^n$, $y>0$, and
\begin{align*}
-H_xu+\frac{1-2\sigma}{y}~u_y+u_{yy}&=0,&\hbox{in }&\Real^n\times(0,\infty);\\
\lim_{y\to0^+}u(x,y)&=f(x),&\hbox{for }&x\in{\mathcal{O}}.
\end{align*}
In addition, for all $x\in{\mathcal{O}}$,
\begin{eqnarray}\label{cosa}
\frac{1}{2\sigma}\lim_{y\to0^+}y^{1-2\sigma}u_y(x,y)=\frac{1}{4^\sigma\Gamma(\sigma)}\int_0^\infty\left(e^{-tH}f(x)-f(x)\right)~\frac{dt}{t^{1+\sigma}}.
\end{eqnarray}
\end{thm}

\begin{proof}
Estimate \eqref{calor LpN} implies that the integral defining $u$ is absolutely convergent and $u_y$ and $u_{yy}$ can be computed by taking the derivatives inside the integral sign. Moreover, by using \eqref{dentro}, we have
$$H_xu(x,y)=\frac{y^{2\sigma}}{4^\sigma\Gamma(\sigma)}\int_0^\infty He^{-tH}f(x)e^{-\frac{y^2}{4t}}~\frac{dt}{t^{1+\sigma}},$$
in the classical sense. Hence, for each $x\in\Real^n$, $u$ verifies the extension problem in the classical sense.
To check that \eqref{cosa} is classical, we begin by recalling that
$$\int_0^\infty e^{-\frac{y^2}{4t}}\left(2\sigma-\frac{y^2}{2t}\right)~\frac{dt}{t^{1+\sigma}}=0.$$
Thus
\begin{align*}
  \frac{1}{2\sigma}~y^{1-2\sigma}u_y(x,y) &= \frac{1}{2\sigma4^\sigma\Gamma(\sigma)}\int_0^\infty e^{-tH}f(x)e^{-\frac{y^2}{4t}}\left(2\sigma-\frac{y^2}{2t}\right)~\frac{dt}{t^{1+\sigma}} \\
   &= \frac{1}{2\sigma4^\sigma\Gamma(\sigma)}\int_0^\infty \left(e^{-tH}f(x)-f(x)\right)e^{-\frac{y^2}{4t}}\left(2\sigma-\frac{y^2}{2t}\right)~\frac{dt}{t^{1+\sigma}}.
\end{align*}
As we shall see later (Remark \ref{util prueba}) the integral in \eqref{cosa} is absolutely convergent for all $x\in{\mathcal{O}}$ and $f$ as in the hypotheses. Therefore \eqref{cosa} follows.
\end{proof}

\begin{rem}\label{util}
In Section \ref{Section Hermite} we will see that for $f\in L^p_N\cap C^2({\mathcal{O}})$, $H^\sigma f$ is well defined and
$$H^\sigma f(x)=\frac{1}{\Gamma(-\sigma)}\int_0^\infty\left(e^{-tH}f(x)-f(x)\right)~\frac{dt}{t^{1+\sigma}},\qquad x\in{\mathcal{O}}.$$
\end{rem}

\begin{rem}
Theorem \ref{Thm:Extension smooth} is valid if $H$ is replaced by $-\Delta$ and the function $f$, with the same smoothness in ${\mathcal{O}}$, belongs to $L_\sigma:=L^1_{n/2+\sigma}$. See the discussion on $(-\Delta)^\sigma$ given in Section \ref{Section Hermite}.
\end{rem}

\begin{lem}[Reflection extension]\label{lem:reflexion}
Fix $R>0$ and $x_0\in\Real^n$. Let $u$ be a solution of
$$-H_xu+\frac{1-2\sigma}{y}~u_y+u_{yy}=0,\qquad\hbox{ in }\Real^n\times(0,R),$$
with
\begin{equation}\label{hypotheses}
\lim_{y\to0^+}y^{1-2\sigma}u_y(x,y)=0,\qquad\hbox{ for every }x\hbox{ such that }\abs{x-x_0}<R.
\end{equation}
Then the extension to $\Real^n\times(-R,R)$ defined by
\begin{equation}\label{eq:reflexion}
\tilde{u}(x,y)=
\left\{
  \begin{array}{ll}
    u(x,y), & y\geq0; \\
    u(x,-y), & y<0;
  \end{array}
\right.
\end{equation}
verifies the degenerate Schr\"{o}dinger equation
\begin{equation}\label{Deg eq}
\dive(\abs{y}^{1-2\sigma}\nabla\tilde{u})-\abs{y}^{1-2\sigma}\abs{x}^2\tilde{u}=0,
\end{equation}
in the weak sense in $B:=\set{(x,y)\in\Real^{n+1}:\abs{x-x_0}^2+y^2<R^2}$.
\end{lem}

\begin{proof}
A nontrivial solution $u$ can be found with the method of Subsection \ref{Subsection Neumann} since the eigenfunctions of the harmonic oscillator $H$ are the Hermite functions $h_\alpha$, $\alpha\in\mathbb{N}_0^n$, with corresponding eigenvalues $\lambda_\alpha=2\abs{\alpha}+n$ (see Section \ref{Section Hermite}). Given $\varphi\in C^\infty_c(B)$ we want to prove that
$$I:=\int_B\left(\nabla\tilde{u}\cdot\nabla\varphi+\abs{x}^2\tilde{u}\varphi\right)\abs{y}^{1-2\sigma}~dx~dy=0.$$
For $\delta>0$ we have
\begin{align*}
  I &= \int_{B\cap\set{\abs{y}\geq\delta}}\dive(\abs{y}^{1-2\sigma}\varphi\nabla\tilde{u})~dx~dy+
\int_{B\cap\set{\abs{y}<\delta}}\left(\nabla\tilde{u}\cdot\nabla\varphi+\abs{x}^2\tilde{u}\varphi\right)\abs{y}^{1-2\sigma}~dx~dy \\
   &= \int_{B\cap\set{\abs{y}=\delta}}\varphi\delta^{1-2\sigma}\tilde{u}_y(x,\delta)~dx+
\int_{B\cap\set{\abs{y}<\delta}}\left(\nabla\tilde{u}\cdot\nabla\varphi+\abs{x}^2\tilde{u}\varphi\right)\abs{y}^{1-2\sigma}~dx~dy.
\end{align*}
As $\delta\to0$, the first term above goes to zero because of \eqref{hypotheses} and the second term goes to zero because $\left(\abs{\nabla\tilde{u}}^2+\abs{x}^2\tilde u\right)\abs{y}^{1-2\sigma}$ is a locally integrable function.
\end{proof}

\begin{proof}[Proof of Theorem \ref{Thm:Harnack}]
Let $u$ be as in Theorem \ref{Thm:Extension smooth}. Since $f$ is a nonnegative function, from \eqref{eq:heat H} and \eqref{Poisson H y Delta} we see that $u\geq0$. Because of Remark \ref{util}, its reflection \eqref{eq:reflexion} satisfies Lemma \ref{lem:reflexion}. Note that \eqref{Deg eq} is a degenerate Schr\"{o}dinger equation with $A_2$ weight $w=\abs{y}^{1-2\sigma}$ and potential $V=\abs{y}^{1-2\sigma}\abs{x}^2$ such that $V/w\in L^p_w$ locally for $p$ large enough. So we can apply the result of \cite{Gutierrez} to obtain the Harnack's inequality for $\tilde{u}$ and thus for $f$.
\end{proof}

\section{Pointwise formula for $H^\sigma$ and some of its consequences}\label{Section Hermite}

The semigroup language adopted in Section \ref{Section Extension}, allows us to get the exact pointwise formula for the fractional Laplacian $(-\Delta)^\sigma$ on $\Real^n.$ The constants involved in the definition are computed exactly in an easy way.

\begin{lem}\label{lem:Frac Lap point}
For $f\in\mathcal{S}$,
\begin{equation}\label{eq:Frac Lap with cnst}
(-\Delta)^\sigma f(x)=\frac{1}{\Gamma(-\sigma)}\int_0^\infty\left(e^{t\Delta}f(x)-f(x)\right)~\frac{dt}{t^{1+\sigma}}= \frac{4^\sigma\Gamma(n/2+\sigma)}{\pi^{n/2}\Gamma(-\sigma)}\PV\int_{\Real^n}\frac{f(z)-f(x)}{\abs{x-z}^{n+2\sigma}}~dz.
\end{equation}
\end{lem}

\begin{proof}
The first identity follows by Fourier transform. From the fact that $e^{t\Delta}1(x)\equiv1$ we can write
\begin{equation}\label{pre Fubini}
\int_0^\infty\left(e^{t\Delta}f(x)-f(x)\right)~\frac{dt}{t^{1+\sigma}}=\int_0^\infty\int_{\Real^n}W_t(x-z)(f(z)-f(x))~dz~\frac{dt}{t^{1+\sigma}}=I_1+I_2,
\end{equation}
where
$$I_1:=\int_0^\infty\int_{\abs{x-z}>1}W_t(x-y)(f(z)-f(x))~dz~\frac{dt}{t^{1+\sigma}},$$
and  $W_t$ is the heat kernel for the Laplacian \eqref{eq:heat Lap}. Use the change of variables $s=\frac{\abs{x-z}^2}{4t}$ to see that
\begin{equation}\label{Frac Lap kernel from heat}
\int_0^\infty\frac{1}{(4\pi t)^{n/2}}~e^{-\frac{\abs{x-z}^2}{4t}}~\frac{dt}{t^{1+\sigma}}=\frac{4^\sigma\Gamma(n/2+\sigma)}{\pi^{n/2}}\cdot\frac{1}{\abs{x-z}^{n+2\sigma}}.
\end{equation}
So, since $f$ is bounded, $I_1$ converges absolutely. Passing to polar coordinates,
$$I_2=\int_0^\infty\frac{1}{(4\pi t)^{n/2}}\int_0^1e^{-\frac{r^2}{4t}}r^{n-1}\int_{\abs{z'}=1}(f(x+rz')-f(x))~dS(z')~dr~\frac{dt}{t^{1+\sigma}}.$$
By Taylor's Theorem, $\int_{\abs{z'}=1}(f(x+rz')-f(x))~dS(z')=C_nr^2\Delta f(x)+O(r^3)$, thus
$$\abs{I_2}\leq C_{n,\Delta f(x)}\int_0^1r^{n+1}\int_0^\infty\frac{e^{-\frac{r^2}{4t}}}{t^{n/2+\sigma}}~\frac{dt}{t}~dr=C_{n,\Delta f(x),\sigma}\int_0^1r^{1-2\sigma}~dr=C_{n,\Delta f(x),\sigma},$$
and $I_2$ converges. Therefore apply Fubini's Theorem in \eqref{pre Fubini} and \eqref{Frac Lap kernel from heat} to get \eqref{eq:Frac Lap with cnst}.
\end{proof}

\begin{rem}
Lemma \ref{lem:Frac Lap point} gives the exact value of the positive constant $c_{n,\sigma}$ in \eqref{Frac Lap Point}. Observe that
\begin{equation}\label{c n sigma}
c_{n,\sigma}=\frac{-4^\sigma\Gamma(n/2+\sigma)}{\pi^{n/2}\Gamma(-\sigma)}\to0,\qquad\hbox{ as }\sigma\to0^+\hbox{ or }\sigma\to1^-.
\end{equation}
\end{rem}

When $f\in\mathcal{S}$ it is clear (by Fourier transform) that $\lim_{\sigma\to1^-}(-\Delta)^\sigma f=-\Delta f$. The next Proposition shows that this is in fact valid for $f\in C^2$. Note that if $f\in\mathcal{S}$ then, from \eqref{Frac Lap Fourier}, $(-\Delta)^\sigma f\notin\mathcal{S}$, but still $(-\Delta)^\sigma f\in C^\infty$. It can be checked that for every $\beta\in\mathbb{N}_0^n$ the function $(1+\abs{x}^{n+2\sigma})D^\beta(-\Delta)^\sigma f(x)$ is bounded. Therefore the set $L_\sigma:=\set{u:\Real^n\to\Real:\norm{u}_{L_\sigma}=\int_{\Real^n}\frac{\abs{u(z)}}{1+\abs{z}^{n+2\sigma}}~dz<\infty}$ (which is $L^1_{n/2+\sigma}$ in \eqref{LpN}), consists of all locally integrable tempered distributions $u$ for which $(-\Delta)^\sigma u$ can be defined. If $f\in L_\sigma$ is $C^2$ in an open set ${\mathcal{O}}$ then it can be proved that $(-\Delta)^\sigma f$ is a continuous function in ${\mathcal{O}}$ and its values are given by the second integral in \eqref{eq:Frac Lap with cnst}. For all the details see \cite{Silvestre Thesis} and \cite{Silvestre CPAM}.

\begin{prop}\label{Prop:Lap sig a 1}
Let $f\in C^2(B_2(x))\cap L^\infty(\Real^n)$ for some $x\in\Real^n$. Then
$$\lim_{\sigma\to1^-}(-\Delta)^\sigma f(x)=-\Delta f(x).$$
\end{prop}

\begin{proof}
Fix an arbitrary $\varepsilon>0$. Since $f\in C^2(B_2(x))$ there exists $\delta=\delta_\varepsilon>0$ such that
\begin{equation}\label{continuidad}
\abs{D^2f(w)-D^2f(w')}<\varepsilon,\qquad\hbox{ for all }w,w'\in\overline{B_1(x)}\hbox{ such that }\abs{w-w'}<\delta.
\end{equation}
Write $(-\Delta)^\sigma f(x)=c_{n,\sigma}(I+II)$ where $I=\int_{\abs{x-z}>\delta}\frac{f(x)-f(z)}{\abs{x-z}^{n+2\sigma}}~dz$. We have $\abs{I}\leq C_n\sigma^{-1}\delta^{-2\sigma}\norm{f}_{L^\infty}$, so that from \eqref{c n sigma}, $c_{n,\sigma}I\to0$ as $\sigma\to1^-$. Using polar coordinates, Taylor's Theorem and recalling that $\int_{\abs{z'}=1}(z_1')^2~dS(z')=\frac{(n/2+1)\pi^{n/2}}{\Gamma(n/2+2)}$,
\begin{align*}
  II &= \int_0^\delta r^{-1-2\sigma}\int_{\abs{z'}=1}\left(f(x)-f(x-rz')\right)~dS(z')~dr \\
   &= \int_0^\delta r^{-1-2\sigma}\int_{\abs{z'}=1}R_1f(x,rz')~dS(z')~dr \\
   &= \int_0^\delta r^{-1-2\sigma}\left[\frac{-\Delta f(x)(n/2+1)\pi^{n/2}r^2}{2\Gamma(n/2+2)}+\int_{\abs{z'}=1} \left(R_1f(x,rz')-\frac{r^2}{2}\langle D^2f(x)z',z'\rangle\right)dS(z')\right]dr \\
   &= \frac{-\Delta f(x)(n/2+1)\pi^{n/2}\delta^{2-2\sigma}}{4\Gamma(n/2+2)(1-\sigma)}+\int_0^\delta r^{-1-2\sigma}\int_{\abs{z'}=1}\left(R_1f(x,rz')-\frac{r^2}{2}\langle D^2f(x)z',z'\rangle\right)dS(z')dr \\ &=:II_1+II_2,
\end{align*}
where $R_1f(x,rz')$ is the Taylor's remainder of first order. Then \eqref{c n sigma} entails
$$c_{n,\sigma}II_1=\frac{-\Delta f(x)\sigma(n/2+1)\Gamma(n/2+\sigma)\delta^{2-2\sigma}}{4^{1-\sigma}\Gamma(n/2+2)\Gamma(2-\sigma)}\to-\Delta f(x)\frac{(n/2+1)\Gamma(n/2+1)}{\Gamma(n/2+2)}=-\Delta f(x),$$
as $\sigma\to1^-$. Finally, by \eqref{continuidad}, $\abs{R_1f(x,rz')-\frac{r^2}{2}\langle D^2f(x)z',z'\rangle}\leq C_nr^2\varepsilon$ and $\abs{II_2}\leq C_n\delta^{2-2\sigma}(1-\sigma)^{-1}\varepsilon$.  Therefore $\lim_{\sigma\to1^{-1}}\abs{c_{n,\sigma}II_2}\leq C_n\varepsilon$.
\end{proof}

\begin{rem}
For $f\in C^2(B_2(x))\cap L_\sigma$ the second identity in \eqref{eq:Frac Lap with cnst} is valid (the idea is to use the continuity of $D^2f$ as in the proof of Proposition \ref{Prop:Lap sig a 1}).
\end{rem}

We shall now discuss the definition of the fractional harmonic oscillator $H^\sigma$ and the pointwise formula for $H^\sigma f(x)$. The eigenfunctions of $H$ (see \cite{Thangavelu}) are the multi-dimensional Hermite functions defined on $\Real^n$ as $h_\alpha(x)=\Phi_\alpha(x)e^{-\abs{x}^2/2}$, $\alpha\in\mathbb{N}_0^n$, where $\Phi_\alpha$ are the multi-dimensional Hermite polynomials, and $Hh_\alpha=(2\abs{\alpha}+n)h_\alpha$. Note that $h_\alpha\in\mathcal{S}$. The set of Hermite functions forms an orthonormal basis of $L^2(\Real^n)$. Let $f\in\mathcal{S}$. The \emph{Hermite series expansion} of $f$ given by
\begin{equation}\label{eq:Series for f}
\sum_\alpha\langle f,h_\alpha\rangle h_\alpha=\sum_{k=0}^\infty\sum_{\abs{\alpha}=k}\langle f,h_\alpha\rangle h_\alpha,
\end{equation}
with $\langle f,h_\alpha\rangle=\int_{\Real^n}fh_\alpha~dx$ (which converges to $f$ in $L^2$), converges uniformly in $\Real^n$ to $f$. This uniform convergence is a consequence of the fact that $\norm{h_\alpha}_{L^\infty(\Real^n)}\leq C$ for all $\alpha\in\mathbb{N}_0^n$ and the following estimate: for every $m\in\mathbb{N}$,
\begin{equation}\label{est coef Herm}
\abs{\langle f,h_\alpha\rangle}=\frac{\abs{\langle H^mf,h_\alpha\rangle}}{(2\abs{\alpha}+n)^m}\leq\frac{\norm{H^mf}_{L^2}}{(2\abs{\alpha}+n)^m},
\end{equation}
since $H$ is a symmetric operator. If $f\in\mathcal{S}$ then
\begin{equation}\label{heat Hermite spectral}
e^{-tH}f(x)=\sum_\alpha e^{-t(2\abs{\alpha}+n)}\langle f,h_\alpha\rangle h_\alpha(x),\qquad t\geq0,
\end{equation}
the series converging uniformly in $\Real^n$. By the given estimates on $\norm{h_\alpha}_{L^\infty}$ and $\abs{\langle f,h_\alpha\rangle}$ the series defining the \emph{fractional Hermite operator}
\begin{equation}\label{eq:Series for Hsigmaf}
H^\sigma f=\sum_\alpha(2\abs{\alpha}+n)^\sigma\langle f,h_\alpha\rangle h_\alpha
\end{equation}
converges uniformly in $\Real^n$.

\begin{lem}\label{lem:H sigma with heat}
For $f\in\mathcal{S}$,
$$H^\sigma f(x)=\frac{1}{\Gamma(-\sigma)}\int_0^\infty\left(e^{-tH}f(x)-f(x)\right)~\frac{dt}{t^{1+\sigma}}.$$
\end{lem}

\begin{proof}
Let $c_\alpha=\langle f,h_\alpha\rangle$. Because of the uniform convergence of the series of \eqref{eq:Series for f}, \eqref{heat Hermite spectral} and \eqref{eq:Series for Hsigmaf} we get
\begin{align*}
  \int_0^\infty\left(e^{-tH}f(x)-f(x)\right)~\frac{dt}{t^{1+\sigma}} &= \int_0^\infty\left(\sum_\alpha e^{-t(2\abs{\alpha}+n)}c_\alpha h_\alpha(x)-\sum_\alpha c_\alpha h_\alpha(x)\right)~\frac{dt}{t^{1+\sigma}} \\
   &= \sum_\alpha c_\alpha h_\alpha(x)\int_0^\infty\left[e^{-t(2\abs{\alpha}+n)}-1\right]~\frac{dt}{t^{1+\sigma}} \\
   &= \Gamma(-\sigma)\sum_\alpha(2\abs{\alpha}+n)^\sigma c_\alpha h_\alpha(x)~=~\Gamma(-\sigma)H^\sigma f(x).
\end{align*}
\end{proof}

We have the following important Lemma whose technical proof is given at the end of this section.

\begin{lem}\label{lem:H sigma distribuciones}
$H^\sigma$ is a continuous operator on $\mathcal{S}$.
\end{lem}

Lemma \ref{lem:H sigma distribuciones} together with the symmetry of $H^\sigma$ on $\mathcal{S}$ (that can be easily verified via Hermite series expansions) allow us to give a distributional definition of $H^\sigma$: for $u\in\mathcal{S}'$, define $H^\sigma u\in\mathcal{S}'$ through
$$\langle H^\sigma u,f\rangle:=\langle u,H^\sigma f\rangle.$$
Therefore $H^\sigma$ is well defined for all functions $u$ that are tempered distributions. In particular, $u$ can be taken from the space $L_N^p$ of \eqref{LpN}, $1\leq p<\infty$, $N>0$.

Recall the expression of $G_t$ given in \eqref{eq:heat H} and the fact that (see \cite{Harboure-deRosa-Segovia-Torrea})
\begin{equation}\label{heat Hermite at 1}
e^{-tH}1(x)=\frac{1}{(\cosh 2t)^{n/2}}~e^{-\frac{\tanh 2t}{2}\abs{x}^2}\leq1.
\end{equation}
Define the nonnegative functions
\begin{equation}\label{def F y B}
F_\sigma(x,z):=\frac{1}{-\Gamma(-\sigma)}\int_0^\infty G_t(x,z)~\frac{dt}{t^{1+\sigma}},\qquad B_\sigma(x):=\frac{1}{\Gamma(-\sigma)}\int_0^\infty\left(e^{-tH}1(x)-1\right)\frac{dt}{t^{1+\sigma}}.
\end{equation}

\begin{thm}\label{Thm:H sig point}
Let $f$ be a function in $L^p_N$ that is $C^2({\mathcal{O}})$ for some open subset ${\mathcal{O}}\subset\Real^n$. Then $H^\sigma f$ is a continuous function in ${\mathcal{O}}$ and
$$H^\sigma f(x)=S_\sigma f(x)+f(x)B_\sigma(x),\qquad x\in{\mathcal{O}},$$
where
\begin{equation}\label{S sigma}
S_\sigma f(x)=\int_{\Real^n}F_\sigma(x,z)(f(x)-f(z))~dz.
\end{equation}
\end{thm}

In \eqref{S sigma} we see that $H^\sigma$ is a nonlocal operator. Before giving the proof of Theorem \ref{Thm:H sig point} we establish some easy consequences.

\begin{thm}[Maximum principle for $H^\sigma$]
Let $f$ be a function in $L^p_N$ that is $C^2$ in an open set ${\mathcal{O}}\subset\Real^n$. Assume that $f\geq0$ and $f(x_0)=0$ for some $x_0\in{\mathcal{O}}$. Then $H^\sigma f(x_0)\leq0$. Moreover, $H^\sigma f(x_0)=0$ only when $f\equiv0$.
\end{thm}

\begin{proof}
By Theorem \ref{Thm:H sig point}, since $f,F_\sigma\geq0$,
$$H^\sigma f(x_0)=\int_{\Real^n}(f(x_0)-f(z))F_\sigma(x_0,z)~dz+f(x_0)B_\sigma(x_0)=-\int_{\Real^n}f(z)F_\sigma(x_0,z)~dz\leq0.$$
If $f(z)>0$ in some set of positive measure, then the last inequality is strict.
\end{proof}

\begin{cor}[Comparison principle for $H^\sigma$]
Let $f,g\in L^p_N\cap C^2({\mathcal{O}})$ be such that $f\geq g$ and $f(x_0)=g(x_0)$ at some $x_0\in{\mathcal{O}}$. Then $H^\sigma f(x_0)\leq H^\sigma g(x_0)$. Moreover, $H^\sigma f(x_0)=H^\sigma g(x_0)$ only when $f\equiv g$.
\end{cor}

We devote the rest of this paper to the proofs of Lemma \ref{lem:H sigma distribuciones}, Theorem \ref{Thm:H sig point}, Proposition \ref{Prop:Heat LpN}, and to complete the missing details at the end of Section \ref{Section Unicidad}.

\begin{proof}[Proof of Lemma \ref{lem:H sigma distribuciones}]
Define the first order partial differential operators
$$A_i:=\frac{\partial}{\partial x_i}+x_i,\quad A_{-i}:=-\frac{\partial}{\partial x_i}+x_i,\qquad i=1,\ldots,n.$$
It is well known that
\begin{equation}\label{Derivadas h alpha}
A_ih_\alpha(x)=(2\alpha_i)^{1/2}h_{\alpha-e_i}(x),\qquad A_{-i}h_\alpha(x)=(2\alpha_i+2)^{1/2}h_{\alpha+e_i}(x),
\end{equation}
where $e_i$ is the $i$th coordinate vector in $\mathbb{N}_0^n$ (see \cite{Thangavelu}). This implies that $H^\sigma f\in C^\infty$ and for all $k\in\mathbb{N}$,
\begin{equation}\label{Derivadas H sigma}
A_{i_1}\cdots A_{i_k}H^\sigma f(x)=\sum_\alpha(2\abs{\alpha}+n)^\sigma\langle f,h_\alpha\rangle A_{i_1}\cdots A_{i_k}h_\alpha(x),\qquad i_l=\pm1,~l=1,\ldots,k,
\end{equation}
the series converging uniformly on $\Real^n$. Since
$$\frac{A_i+A_{-i}}{2}=x_i,\quad\frac{A_i-A_{-i}}{2}=\frac{\partial}{\partial x_i},\qquad i=1,\ldots,n,$$
for each multi-index $\gamma,\beta\in\mathbb{N}_0^n$ we can write $x^\gamma D^\beta=x_1^{\gamma_1}\cdots x_n^{\gamma_n}\frac{\partial^{\abs{\beta}}}{\partial x_1^{\beta_1}\cdots\partial x_n^{\beta_n}}$ as a finite linear combination of operators $A_i$ and $A_{-i}$. Therefore, to check that $x^\gamma D^\beta H^\sigma f\in L^\infty$, it is enough to verify that for each $k\in\mathbb{N}$, $A_{i_1}\cdots A_{i_k}H^\sigma f\in L^\infty$ where $\set{i_1,\ldots,i_k}\subset\set{-1,1}$. The identities in \eqref{Derivadas h alpha} easily imply the following commutation relations for Hermite functions and thus for $f\in S$:
$$
\left\{
  \begin{array}{ll}
    A_iH^\sigma f=(H+2)^\sigma A_if, & 1\leq i\leq n; \\
    A_iH^\sigma f=(H-2)^\sigma A_if, & -n\leq i\leq-1.
   \end{array}
\right.
$$
Here $(H\pm2)^{\sigma}A_if:=\sum_\alpha(2\abs{\alpha}+n\pm2)^\sigma\langle A_if,h_\alpha\rangle h_\alpha$. Hence, in \eqref{Derivadas H sigma},
$$A_{i_1}\cdots A_{i_k}H^\sigma f=\sum_\alpha(2\abs{\alpha}+n+2j)^\sigma\langle g,h_\alpha\rangle h_\alpha,$$
for some $j\in\mathbb{Z}$ and $g:=A_{i_1}\cdots A_{i_k}f\in\mathcal{S}$. For $m\in\mathbb{N}$ sufficiently large, as in \eqref{est coef Herm}, we have
$$\abs{\sum_\alpha(2\abs{\alpha}+n+2j)^\sigma\langle g,h_\alpha\rangle h_\alpha(x)}\leq\norm{H^mg}_{L^2(\Real^n)}\sum_\alpha\frac{(2\abs{\alpha}+n+2j)^\sigma}{(2\abs{\alpha}+n)^m}=C\norm{H^mg}_{L^2(\Real^n)}.$$
Therefore $x^\gamma D^\beta H^\sigma f\in L^\infty$. Moreover,
\begin{align*}
  \abs{x^\gamma D^\beta H^\sigma f(x)} &= \abs{\sum c_{i,k}A_{i_1}\cdots A_{i_k}H^\sigma f(x)}~\leq~C\sum\abs{(H+2j)^\sigma A_{i_1}\cdots A_{i_k}f(x)} \\
   &\leq C\left(\hbox{seminorms in }\mathcal{S}\hbox{ of }(A_{i_1}\cdots A_{i_k}f)\right)=C\left(\hbox{seminorms in }\mathcal{S}\hbox{ of }f\right).
\end{align*}
\end{proof}

For the proof of Theorem \ref{Thm:H sig point} we need some estimates on $G_t$, $F_\sigma$ and $B_\sigma$. First we derive some equivalent formulas for these kernels. Consider the change of parameters due to S. Meda
\begin{equation}\label{Meda transf}
t=t(s)=\frac{1}{2}\log\frac{1+s}{1-s},\qquad t\in(0,\infty),~s\in(0,1),
\end{equation}
that produces
\begin{equation}\label{def dmu sigma}
\frac{dt}{t^{1+\sigma}}=d\mu_\sigma(s):=\frac{ds}{(1-s^2)\left(\frac{1}{2}\log\frac{1+s}{1-s}\right)^{1+\sigma}},\qquad t\in(0,\infty),~s\in(0,1).
\end{equation}
Then the heat kernel in \eqref{eq:heat H} can be written as
$$G_{t(s)}(x,z)=\left(\frac{1-s^2}{4\pi s}\right)^{n/2}e^{-\frac{1}{4}\left[s\abs{x+z}^2+\frac{1}{s}\abs{x-z}^2\right]},\quad s\in(0,1),$$
and, from \eqref{def F y B} and \eqref{def dmu sigma},
$$F_\sigma(x,z)=\frac{1}{-\Gamma(-\sigma)}\int_0^1 G_{t(s)}(x,z)~d\mu_\sigma(s),\quad B_\sigma(x)=\frac{1}{\Gamma(-\sigma)}\int_0^1\left(e^{-t(s)H}1(x)-1\right)d\mu_\sigma(s).$$

\begin{lem}\label{lem:Heat krn est}
For all $s\in(0,1)$ and $x,z\in\Real^n$,
\begin{equation}\label{Heat krnl estm}
G_{t(s)}(x,z)\leq C\left(\frac{1-s}{s}\right)^{n/2}e^{-\frac{\abs{x}\abs{x-z}}{C}}e^{-\frac{\abs{x-z}^2}{Cs}}.
\end{equation}
In particular,
\begin{equation}\label{Heat krnl estm2}
G_{t(s)}(x,z)\leq\frac{C}{\abs{x-z}^n}~(1-s)^{n/2}e^{-\frac{\abs{x}\abs{x-z}}{C}}e^{-\frac{\abs{x-z}^2}{C}}e^{-\frac{\abs{x-z}^2}{Cs}}.
\end{equation}
\end{lem}

\begin{proof}
The second estimate in the statement follows immediately from \eqref{Heat krnl estm}. Note that
$$G_{t(s)}(x,z)\leq C\left(\frac{1-s}{s}\right)^{n/2}e^{-\frac{\abs{x-z}^2}{8s}}e^{-\frac{1}{8}\left[s\abs{x+z}^2+\frac{1}{s}\abs{x-z}^2\right]}\leq C\left(\frac{1-s}{s}\right)^{n/2}e^{-\frac{\abs{x-z}^2}{8s}}e^{-\frac{1}{8}\abs{x-z}\abs{x+z}}.$$
We prove the second inequality above. Assume first that $\abs{x-z}\leq\abs{x+z}$. Then by minimizing the function $\theta(s):=\frac{s}{8}\abs{x+z}^2+\frac{1}{8s}\abs{x-z}^2$ for $s\in(0,1)$ we get
$e^{-\frac{1}{8}\left[s\abs{x+z}^2+\frac{1}{s}\abs{x-z}^2\right]}\leq e^{-\frac{1}{8}\abs{x-z}\abs{x+z}}$. In the case $\abs{x+z}<\abs{x-z}$ we have $e^{-\frac{1}{8}\left[s\abs{x+z}^2+\frac{1}{s}\abs{x-z}^2\right]}\leq e^{-\frac{1}{8s}\abs{x-z}^2}=e^{-\frac{1}{8s}\abs{x-z}\abs{x-z}}\leq e^{-\frac{1}{8}\abs{x-z}\abs{x+z}}$, for all $s\in(0,1)$. To obtain estimate \eqref{Heat krnl estm} proceed as follows: if $x\cdot z>0$ then $\abs{x+z}\geq\abs{x}$ which gives $e^{-\frac{1}{8}\abs{x-z}\abs{x+z}}\leq e^{-\frac{1}{8}\abs{x}\abs{x-z}}$; if $x\cdot z\leq0$ then $\abs{x-z}\geq\abs{x}$ and in this situation
$$e^{-\frac{\abs{x-z}^2}{8s}}e^{-\frac{1}{8}\abs{x-z}\abs{x+z}}\leq e^{-\frac{\abs{x-z}^2}{16s}}e^{-\frac{\abs{x}\abs{x-z}}{16s}}\leq e^{-\frac{\abs{x-z}^2}{16s}}e^{-\frac{\abs{x}\abs{x-z}}{16}}.$$
\end{proof}

Observe in \eqref{def dmu sigma} that
\begin{equation}\label{est dmu rho}
d\mu_\sigma(s)\sim\frac{ds}{s^{1+\sigma}},~s\sim 0,\qquad d\mu_\sigma(s)\sim\frac{ds}{(1-s)(-\log(1-s))^{1+\sigma}},~s\sim1.
\end{equation}

\begin{lem}\label{lem:F and B est}
For all $x,z\in\Real^n$,
\begin{equation}\label{F sigma estm}
F_\sigma(x,z)\leq\frac{C}{\abs{x-z}^{n+2\sigma}}~e^{-\frac{\abs{x}\abs{x-z}}{C}}e^{-\frac{\abs{x-z}^2}{C}}\quad\hbox{ and }\quad B_\sigma(x)\leq C\left(1+\abs{x}^{2\sigma}\right).
\end{equation}
Moreover, $B_\sigma\in C^\infty(\Real^n)$.
\end{lem}

\begin{proof}
Estimate \eqref{Heat krnl estm2} gives
$$F_\sigma(x,z)\leq C\frac{e^{-\frac{\abs{x}\abs{x-z}}{C}}}{\abs{x-z}^n}\int_0^1(1-s)^{n/2}e^{-\frac{\abs{x-z}^2}{Cs}}~d\mu_\sigma(s).$$
Then \eqref{est dmu rho} implies
$$
\int_0^{1/2}(1-s)^{n/2}e^{-\frac{\abs{x-z}^2}{Cs}}~d\mu_\sigma(s)\leq C\int_0^{1/2}e^{-\frac{\abs{x-z}^2}{Cs}}\frac{ds}{s^{1+\sigma}}\leq C
\left\{
    \begin{array}{ll}
        \frac{1}{\abs{x-z}^{2\sigma}}, & \hbox{ if }\abs{x-z}<1; \\
        e^{-\frac{\abs{x-z}^2}{C}}, & \hbox{ if }\abs{x-z}\geq1;
    \end{array}
\right.
$$
and
$$\int_{1/2}^1(1-s)^{n/2}e^{-\frac{\abs{x-z}^2}{Cs}}~d\mu_\sigma(s)\leq e^{-\frac{\abs{x-z}^2}{C}}\int_{1/2}^1\frac{ds}{(1-s)(-\log(1-s))^{1+\sigma}}=Ce^{-\frac{\abs{x-z}^2}{C}},$$
thus the first inequality in \eqref{F sigma estm} follows.

Apply \eqref{Meda transf} in \eqref{heat Hermite at 1} to obtain
\begin{equation}\label{heat Hermite at 1 con s}
e^{-t(s)H}1(x)=\left(\frac{1-s^2}{1+s^2}\right)^{n/2}e^{-\frac{s}{1+s^2}\abs{x}^2}.
\end{equation}
Then, up to the factor $\frac{1}{-\Gamma(-\sigma)}$, we can write
$$B_\sigma(x)=\int_0^1\left[\left(\frac{1-s^2}{1+s^2}\right)^{n/2}-1\right]e^{-\frac{s}{1+s^2}\abs{x}^2}~d\mu_\sigma(s)+\int_0^1\left(e^{-\frac{s}{1+s^2}\abs{x}^2}-1\right)~d\mu_\sigma(s)=I+II.$$
To estimate $I$ and $II$ we use \eqref{est dmu rho} and the Mean Value Theorem. That is,
$$\abs{I}\leq C\int_0^{1/2}\abs{\left(\frac{1-s^2}{1+s^2}\right)^{n/2}-1}\frac{ds}{s^{\sigma+1}}+\int_{1/2}^1~d\mu_\sigma(s)\leq C\int_0^{1/2}s^2~\frac{ds}{s^{1+\sigma}}+C~=~C.$$
For $II$ we consider two cases. Assume first that $\abs{x}^2\leq2$. Then
$$\abs{II}\leq C\int_0^{1/2}\abs{e^{-\frac{s}{1+s^2}\abs{x}^2}-1}~\frac{ds}{s^{1+\sigma}}+\int_{1/2}^1~d\mu_\sigma(s)\leq C\int_0^{1/2}\abs{x}^2s~\frac{ds}{s^{1+\sigma}}+C\leq C.$$
In the case $\abs{x}^2>2$,
\begin{align*}
  \abs{II} &\leq \abs{x}^2\int_0^{\frac{1}{\abs{x}^2}}s~\frac{ds}{s^{1+\sigma}}+\int_{\frac{1}{\abs{x}^2}}^1~d\mu_\sigma(s)\leq \abs{x}^2\int_0^{\frac{1}{\abs{x}^2}}s^{-\sigma}~ds+\int_\frac{1}{\abs{x}^2}^1\frac{ds}{(1-s)\left(-\log(1-s)\right)^{1+\sigma}} \\
   &= C\abs{x}^{2\sigma}+C\left[-\log\left(1-\frac{1}{\abs{x}^2}\right)\right]^{-\sigma}\leq C\abs{x}^{2\sigma},
\end{align*}
since $-\log(1-s)\sim s$ as $s\to0$. Therefore the second fact of \eqref{F sigma estm} follows. $B_\sigma$ is differentiable since the gradient of the integrand in its definition is bounded by $$2\abs{x}\frac{s}{1+s^2}\left(\frac{1-s^2}{1+s^2}\right)^{n/2}e^{-\frac{s}{1+s^2}\abs{x}^2}\leq C\abs{x}s\in L^1((0,1);d\mu_\sigma(s)),$$
thus we can differentiate inside the integral:
$$\nabla B_\sigma(x)=2x\int_0^1\frac{s}{1+s^2}\left(\frac{1-s^2}{1+s^2}\right)^{n/2}e^{-\frac{s}{1+s^2}\abs{x}^2}~d\mu_\sigma(s).$$
For higher order derivatives we can proceed similarly.
\end{proof}

\begin{proof}[Proof of Theorem \ref{Thm:H sig point}]
Take first $f\in\mathcal{S}$. Since $e^{-tH}1(x)$ is not a constant function we write
\begin{align*}
  \lefteqn{\int_0^\infty\left(e^{-tH}f(x)-f(x)\right)~\frac{dt}{t^{1+\sigma}}=\int_0^\infty \left(\int_{\Real^n}G_t(x,z)f(z)~dz-f(x)\right)~\frac{dt}{t^{1+\sigma}}} \\
   &= \int_0^\infty\left[\int_{\Real^n}G_t(x,z)(f(z)-f(x))~dz+f(x)\left(\int_{\Real^n} G_t(x,z)~dz-1\right)\right]~\frac{dt}{t^{1+\sigma}} \\
   &= \int_0^\infty\int_{\Real^n}G_t(x,z)(f(z)-f(x))~dz~\frac{dt}{t^{1+\sigma}}+ f(x)\int_0^\infty\left(e^{-tH}1(x)-1\right)\frac{dt}{t^{1+\sigma}} \\
   &= \int_0^1\int_{\Real^n}G_{t(s)}(x,z)(f(z)-f(x))~dz~d\mu_\sigma(s)+f(x)B_\sigma(x).
\end{align*}
Due to Lemma \ref{lem:H sigma with heat}, the first integral above is well defined and converges absolutely. Write the integral in the last line as $I_\delta+I_{\delta^c}$ with $I_{\delta^c}=\int_0^1\int_{\abs{x-z}>\delta}G_{t(s)}(x,z)(f(z)-f(x))~dz~d\mu_\sigma(s)$, for some $\delta>0$ (in this step $\delta$ is arbitrary, but we will fix it later). Estimate \eqref{F sigma estm} implies that $I_{\delta^c}$ is absolutely convergent and $\abs{I_{\delta^c}}\leq C\norm{f}_{L^\infty(\Real^n)}$. Pass to polar coordinates in $I_\delta$:
\begin{align*}
  I_\delta &= \int_0^1\left(\frac{1-s^2}{4\pi s}\right)^{n/2}\int_{\abs{x-z}<\delta}e^{-\frac{1}{4}\left[s\abs{x+z}^2+\frac{1}{s}\abs{x-z}^2\right]} (f(z)-f(x))~dz~d\mu_\sigma(s) \\
   &= \int_0^1\left(\frac{1-s^2}{4\pi s}\right)^{n/2}\int_0^\delta r^{n-1}e^{-\frac{r^2}{4s}}\int_{\abs{z'}=1}e^{-\frac{s}{4}\abs{2x+rz'}^2}(f(x+rz')-f(x))~dS(z')drd\mu_\sigma(s).
\end{align*}
To estimate $I_{S^{n-1}}:=\int_{\abs{z'}=1}e^{-\frac{s}{4}\abs{2x+rz'}^2}(f(x+rz')-f(x))~dS(z')$ use the Taylor expansions of $f$ and $\psi_s(w):=e^{-\frac{s}{4}\abs{w}^2}$ and cancel out terms:
\begin{align*}
  I_{S^{n-1}} &= \int_{\abs{z'}=1}\left(e^{-\frac{s}{4}\abs{2x}^2}+R_0\psi_s(x,rz')\right)\left(\nabla f(x)(rz')+R_1f(x,rz')\right)~dS(z') \\
   &= \int_{\abs{z'}=1}\left[e^{-\frac{s}{4}\abs{2x}^2}R_1f(x,rz')+R_0\psi_s(x,rz')\nabla f(x)(rz')+R_0\psi_s(x,rz')R_1f(x,rz')\right]dS(z').
\end{align*}
Since $\abs{R_0\psi_s(x,rz')}\leq s^{1/2}r$ and $\abs{R_1f(x,rz')}\leq\norm{D^2f}_{L^\infty(B_\delta(x))}r^2$, we have $\abs{I_{S^{n-1}}}\leq Cr^2$. Thus
\begin{align*}
  \abs{I_\delta} &\leq \int_0^1\left(\frac{1-s^2}{4\pi s}\right)^{n/2}\int_0^\delta r^{n-1}e^{-\frac{r^2}{4s}}\abs{I_{S^{n-1}}}~dr~d\mu_\sigma(s) \\
   &\leq C\int_0^\delta r^{n+1}\int_0^1\frac{1}{s^{n/2}}~e^{-\frac{r^2}{4s}}~d\mu_\sigma(s)~dr\leq C\int_0^\delta r^{n+1}\frac{1}{r^{n+2\sigma}}~dr=C\delta^{2-2\sigma}.
\end{align*}
Hence $I_\delta$ converges. The conclusion follows, for $f\in\mathcal{S}$, by Fubini's Theorem.

Now assume that $f\in L^p_N$, $1\leq p<\infty$, $N>0$, is a $C^2$ function in ${\mathcal{O}}$. Then $H^\sigma f$ is well defined as a tempered distribution. Fix an arbitrary $x\in{\mathcal{O}}$ and take $\delta>0$ so that $B_\delta(x)\subset{\mathcal{O}}$. Observe that the integral in \eqref{S sigma} is well defined: just apply Taylor's Theorem (as above) in $I_{\delta}$ and the $L^p_N$ condition together with \eqref{F sigma estm} in $I_{\delta^c}$. Let $\varepsilon>0$. There exists a sequence $f_k\in C^\infty_c(\Real^n)$ such that $\norm{D^2f_k}_{L^\infty(B_\delta(x))}\leq\norm{D^2f}_{L^\infty(B_\delta(x))}$ for all $k$, $f_k$ converges uniformly to $f$ in $B_\delta(x)$ and $f_k\to f$ in the norm of $L^p_N$, as $k\to\infty$ (use mollifiers and multiplication by a smooth cutoff function). Since $B_\sigma$ is a continuous function, $f_kB_\sigma$ converges uniformly to $fB_\sigma$ on $B_\delta(x)$. Let $0<\rho<\delta/2$ be such that for all $k$
$$\abs{\int_{B_\rho(x)}F_\sigma(x,z)(f_k(x)-f_k(z))~dz}<\frac{\varepsilon}{3},\quad\hbox{ and }\quad\abs{\int_{B_\rho(x)}F_\sigma(x,z)(f(x)-f(z))~dz}<\frac{\varepsilon}{3}.$$
For $k$ sufficiently large, by H\"{o}lder's inequality,
\begin{align*}
  \lefteqn{\abs{\int_{B_\rho^c(x)}F_\sigma(x,z)(f_k(x)-f_k(z))~dz-\int_{B_\rho^c(x)}F_\sigma(x,z)(f(x)-f(z))~dy}\leq} \\
   &\leq \abs{f_k(x)-f(x)}\int_{B_\rho^c(x)}F_\sigma(x,z)~dz+\int_{B_\rho^c(x)}F_\sigma(x,z)\abs{f_k(z)-f(z)}~dz \\
   &\leq C\left(\abs{f_k(x)-f(x)}+\norm{f_k-f}_{L^p_N}\right)<\frac{\varepsilon}{3}.
\end{align*}
Thus
$$S_\sigma f_k(x)\rightrightarrows\int_{\Real^n}F_\sigma(x,z)(f(x)-f(z))~dz$$
in $B_\delta(x)$. But $H^\sigma f_k\to H^\sigma f$ in $\mathcal{S}'$. By uniqueness of the limits, $S_\sigma f(x)$ coincides with the integral in \eqref{S sigma}. Moreover, $H^\sigma f$ is continuous in $B_\delta(x)$ because it is the uniform limit of continuous functions.
\end{proof}

\begin{proof}[Proof of Proposition \ref{Prop:Heat LpN}]
By \eqref{Meda transf}, \eqref{Heat krnl estm} and H\"older's inequality,
$$\abs{e^{-t(s)H}f(x)}\leq \frac{C\norm{f}_{L^p_N}}{s^{n/2}}\left(\int_{\Real^n}e^{-\frac{p'\abs{x-z}^2}{C}}(1+\abs{z}^2)^{Np'}~dz\right)^{1/p'}\leq C\frac{(1+\abs{x}^\rho)\norm{f}_{L^p_N}}{s^{n/2}}.$$
For \eqref{calor LpN} note that if $0<s<\frac{1}{2}$, then $s<t(s)<\frac{4}{3}s$. The equality $\partial_te^{-tH}f(x)=\int_{\Real^n}\partial_tG_t(x,z)f(z)~dz$ is valid if the last integral is absolutely convergent for all $t$ in some interval. But $\partial_tG_t(x,z)f(z)=-H_xG_t(x,z)f(z)$, therefore we have to verify that $\int_{\Real^n}H_xG_{t(s)}(x,z)f(z)~dz$ converges absolutely for all $s$ in some interval. This last statement is true since
$$\abs{\nabla_xG_{t(s)}(x,z)}\leq\left(\frac{1-s^2}{s}\right)^{n/2}\frac{1}{s^{1/2}}~e^{-c\left[s\abs{x+z}^2+\frac{1}{s}\abs{x-z}^2\right]}$$
and
$$\abs{D_x^2G_{t(s)}(x,z)}\leq\left(\frac{1-s^2}{s}\right)^{n/2}\frac{1}{s}~e^{-c\left[s\abs{x+z}^2+\frac{1}{s}\abs{x-z}^2\right]},$$
which give estimates similar to \eqref{Heat krnl estm} for $\nabla G_{t(s)}$ and $D^2G_{t(s)}$. Hence $\partial_te^{-tH}f(x)=-H_xe^{-tH}f(x)$ and \eqref{dentro} follows. Observe that $t(s)\to0$ if and only if $s\to0$. For $x\in{\mathcal{O}}$ we have
$$\abs{e^{-t(s)H}f(x)-f(x)}\leq\abs{\int_{\Real^n}G_{t(s)}(x,z)(f(z)-f(x))~dz}+\abs{f(x)}\abs{e^{-t(s)H}1(x)-1}.$$
The last term above tends to $0$ as $t(s)\to0$ because of \eqref{heat Hermite at 1 con s}. Let $\delta>0$ be such that $B_\delta(x)\subset{\mathcal{O}}$. Then, as $f\in C^1(\overline{B_\delta(x)})$,
$$\abs{\int_{B_\delta(x)}G_{t(s)}(x,z)(f(z)-f(x))~dz}\leq C\int_{B_\delta(x)}\frac{e^{-\frac{\abs{x-y}^2}{Cs}}}{s^{n/2}}\abs{z-x}dz\leq C\int_{B_\delta(x)}\frac{e^{-\frac{\abs{x-y}^2}{Cs}}}{\abs{z-x}^{n-1}}~dz\to0,$$
when $s\to0$, by the Dominated Convergence Theorem. On the other hand,
\begin{align*}
  \abs{\int_{B_\delta^c(x)}G_{t(s)}(x,z)(f(z)-f(x))~dz} &\leq \left(\int_{B_\delta^c(x)}e^{-\frac{p\abs{x-z}^2}{2C}}\left(\frac{\abs{f(z)}^p}{(1+\abs{z}^2)^{Np}} +\frac{\abs{f(x)}^p}{(1+\abs{z}^2)^{Np}}\right)dz\right)^{1/p} \\
  & \qquad\times \frac{C}{s^{n/2}}\left(\int_{B_\delta^c(x)}e^{-\frac{p'\abs{x-z}^2}{Cs}}e^{-\frac{p'\abs{x-z}^2}{2C}} (1+\abs{z}^2)^{Np'}~dz\right)^{1/p'} \\
  &=: I\times II.
\end{align*}
Clearly $I<\infty$ and, by dominated convergence,
$$II\leq C\left(\int_{B_\delta^c(x)}\frac{e^{-\frac{p'\abs{x-z}^2}{Cs}}}{\abs{x-z}^{np'}}~e^{-\frac{p'\abs{x-z}^2}{C}}(1+\abs{z}^2)^{Np'}~dz\right)^{1/p'}\to0,\qquad\hbox{ as }s\to0.$$
\end{proof}

\begin{rem}\label{util prueba}
If $f\in L^p_N\cap C^2({\mathcal{O}})$ then, for each $x\in{\mathcal{O}}$,
$$\int_0^\infty\abs{e^{-tH}f(x)-f(x)}~\frac{dt}{t^{1+\sigma}}=\int_0^1\abs{e^{-t(s)H}f(x)-f(x)}~d\mu_\sigma(s)<\infty.$$
Indeed, by \eqref{calor LpN},
$$\int_{\frac{1}{2}\log3}^\infty\abs{e^{-tH}f(x)-f(x)}~\frac{dt}{t^{1+\sigma}}\leq C(x)\int_{\frac{1}{2}\log3}^\infty\frac{1}{t^{n/2}}~\frac{dt}{t^{1+\sigma}}<\infty,$$
and
\begin{align*}
  \lefteqn{\int_0^{\frac{1}{2}\log3}\abs{e^{-tH}f(x)-f(x)}~\frac{dt}{t^{1+\sigma}} =\int_0^{1/2}\abs{e^{-t(s)H}f(x)-f(x)}~d\mu_\sigma(s)} \\
   &\leq C\int_0^{1/2}\abs{\int_{\Real^n}G_{t(s)}(x,z)(f(z)-f(x))dz}\frac{ds}{s^{1+\sigma}}+ C\abs{f(x)}\int_0^{1/2}\left[1-\left(\frac{1-s^2}{1+s^2}\right)^{n/2}e^{-\frac{s}{1+s^2}\abs{x}}\right] \frac{ds}{s^{1+\sigma}}.
\end{align*}
Both integrals above are finite: the first one by the arguments in the proof of Theorem \ref{Thm:H sig point} (Taylor's Theorem) and the second one because of the Mean Value Theorem.
\end{rem}

\noindent\textbf{Acknowledgments.} We are very grateful to the referee for his detailed comments. The variety of his substantial suggestions certainly helped us to improve the results and presentation of the paper in an essential way.



\end{document}